\documentclass{amsart}
\usepackage{enumerate, amssymb}

\usepackage{diagrams}
\diagramstyle[scriptlabels, small]

\def\pdfsyncstart{}
\def\pdfsyncstop{}
\def\bdi{\pdfsyncstop\begin{diagram}}
\def\edi{\end{diagram}\pdfsyncstart}

\newarrow{Equal}{}{=}{}{=}{}

\def\lto{{\longrightarrow}}
\def\hto{{\hookrightarrow}}
\def\xto{\xrightarrow}

\newcommand{\fX}{{\mathfrak X}}

\newcommand{\fZ}{{\mathfrak Z}}

\newcommand{\cala}{{\mathcal A}}
\newcommand{\calb}{{\mathcal B}}
\newcommand{\calc}{{\mathcal C}}

\newcommand{\cale}{{\mathcal E}}
\newcommand{\calf}{{\mathcal F}}
\newcommand{\calg}{{\mathcal G}}

\newcommand{\cali}{{\mathcal I}}
\newcommand{\calj}{{\mathcal J}}

\newcommand{\calm}{{\mathcal M}}
\newcommand{\caln}{{\mathcal N}}
\newcommand{\calo}{{\mathcal O}}
\newcommand{\calp}{{\mathcal P}}
\newcommand{\calq}{{\mathcal Q}}
\newcommand{\calr}{{\mathcal R}}
\newcommand{\cals}{{\mathcal S}}
\newcommand{\calt}{{\mathcal T}}

\newcommand{\calov}{{\calo_V}}

\newcommand{\calox}{{\calo_X}}
\newcommand{\caloy}{{\calo_Y}}

\def\calas{{\cala_{*}}}
\def\calbs{{\calb_{*}}}    

\def\cales{{\cale_{*}}}
\def\calfs{{\calf_{*}}}
\def\calgs{{\calg_{*}}}

\def\calms{{\calm_{*}}}
\def\calns{{\caln_{*}}}
\def\calps{{\calp_{*}}}

\def\calrs{{\calr_{*}}}
\def\calss{{\cals_{*}}}

\def\caloxs{{\calo_{X_{*}}}}
\def\caloys{{\calo_{Y_{*}}}}
\def\calozs{{\calo_{Z_{*}}}}
\def\calows{{\calo_{W_{*}}}}
\def\Omegas{{\Omega_{*}}}

\newcommand{\bbbh}{{\mathbb H}}
\newcommand{\bbbl}{{\mathbb L}}
\newcommand{\bbbq}{{\mathbb Q}}
\newcommand{\bbbc}{{\mathbb C}}

\newcommand{\bbbs}{{\mathbb S}}
\newcommand{\bbbn}{{\mathbb N}}

\newcommand{\bbbv}{{\mathbb V}}
\newcommand{\bbbz}{{\mathbb Z}}

\newcommand{\bdot}{{\scriptscriptstyle\bullet}}

\newcommand{\vp}{\varphi}

\newcommand{\lr}{\,\left\lrcorner\right.\,}

\DeclareMathOperator{\At}{At}

\DeclareMathOperator{\ch}{ch}

\DeclareMathOperator{\dotimes}{\underline{\otimes}}

\DeclareMathOperator{\Ext}{Ext}

\DeclareMathOperator{\HH}{HH}

\DeclareMathOperator{\Hom}{Hom}
\DeclareMathOperator{\id}{id}

\DeclareMathOperator{\Ker}{Ker}
\DeclareMathOperator*{\liminj}{\underrightarrow{\lim\vphantom{(}}}

\DeclareMathOperator{\Spec}{Spec}

\DeclareMathOperator{\Tor}{Tor}

\DeclareMathOperator{\Tr}{Tr}
\DeclareMathOperator{\ex}{ex}

\DeclareMathOperator{\AH}{AH}
\DeclareMathOperator{\AC}{AC}

\def\nablau{{\nabla_u}} 
\def\Atu{{\At_u}} 

\def\tAt{\widetilde{\At}}

\newcommand{{\sbullet}}{{\scriptstyle\bullet}}

\def\ts{{\tilde\bbbs}}

\def\Cech{\v{C}ech\ }

\def\calcb{{\calc^{\bdot}}}

\theoremstyle{definition}
\newtheorem{defn}{Definition}[subsection]
\newtheorem{con}[defn]{Convention}

\theoremstyle{plain}
\newtheorem{prop}[defn]{Proposition}
\newtheorem{theorem}[defn]{Theorem}
\newtheorem{lem}[defn]{Lemma}
\newtheorem{cor}[defn]{Corollary}

\theoremstyle{remark}
\newtheorem{rem}[defn]{Remark}

\newtheorem{sit}[defn]{}

\def\bpro{\begin{prop}}
\def\epro{\end{prop}}

\def\bcor{\begin{cor}}
\def\ecor{\end{cor}}

\def\bthm{\begin{theorem}}
\def\ethm{\end{theorem}}

\def\bdfn{\begin{defn}}
\def\edfn{\end{defn}}

\def\brem{\begin{rem}}
\def\erem{\end{rem}}

\def\bsit{\begin{sit}}
\def\esit{\end{sit}}

\def\blem{\begin{lem}}
\def\elem{\end{lem}}

\def\bdi{\pdfsyncstop\begin{diagram}}
\def\edi{\end{diagram}\pdfsyncstart}

\def\ba{\begin{array}}
\def\ea{\end{array}}

\def\bnum{\begin{enumerate}}
\def\enum{\end{enumerate}}

\def\be{\begin{equation}}
\def\ee{\end{equation}}

\def\bproof{\begin{proof}}
\def\eproof{\end{proof}}

\begin{document}
\title[Decomposition of Hochschild (Co-)Homology]{The global
decomposition theorem for Hochschild (co-)homology of singular spaces
via the Atiyah--Chern character}

\author{Ragnar-Olaf Buchweitz}
\address{Dept.\ of Computer and Mathematical Sciences, University of
Toronto at Scarborough, Toronto, Ont.\ M1C 1A4, Canada}
\email{ragnar@math.toronto.edu}

\author{Hubert Flenner}
\address{Fakult\"at f\"ur Mathematik der Ruhr-Universit\"at,
Universit\"atsstr.\ 150, Geb.\ NA 2/72, 44780 Bochum, Germany}
\email{Hubert.Flenner@rub.de}

\thanks{The authors were partly supported by NSERC grant
3-642-114-80, by the DFG Schwerpunkt ``Global Methods in Complex
Geometry", and by the Research in Pairs Program at
Math.~Forschungsinstitut Oberwolfach funded by the
Volkswagen Foundation.}

\begin{abstract} 
We generalize the decomposition theorem of Hochschild, Kostant and 
Rosenberg for Hochschild (co-)homology to arbitrary morphisms between 
complex spaces or schemes over a field of characteristic zero. 
To be precise, we show that for each such morphism $X\to Y$, 
the Hochschild complex $\bbbh_{X/Y}$, as introduced in \cite{BFl2}, 
decomposes naturally in the derived category $D(X)$ into 
$\bigoplus_{p\ge 0}\bbbs^p(\bbbl_{X/Y}[1])$, the direct sum of the 
derived symmetric powers of the shifted cotangent complex, a result
due to Quillen in the affine case. 

Even in the affine case, our proof is new and provides further information.
It shows that the decomposition is given explicitly and naturally by the 
{\em universal Atiyah--Chern character\/}, the exponential of the universal Atiyah class.

We further use the decomposition theorem to show that the semiregularity map
for perfect complexes factors through Hochschild homology and, in turn, factors
the {\em Atiyah--Hochschild character\/} through the characteristic homomorphism from
Hochschild cohomology to the graded centre of the derived category.
\end{abstract}

\maketitle

{\footnotesize\tableofcontents}

\section*{Introduction}

This paper builds on \cite{BFl2}, where we introduced a 
Hochschild complex $\bbbh_{X/Y}$ for any morphism of analytic spaces or 
arbitrary schemes $f:X\to Y$ and defined Hochschild (co-)homology  for a complex 
$\calm$ in $D(X)$ respectively as
\begin{align*}
\HH_{\bdot}^{X/Y}(\calm)&:=
\Tor_{\bdot}^X(\bbbh_{X/Y}, \calm)\quad\text{and}\quad
\HH^{\bdot}_{X/Y}(\calm):=
\Ext^{\bdot}_X(\bbbh_{X/Y}, \calm)\,.
\end{align*}
It also extends considerably our earlier work \cite{BFl} in that it broadens 
vastly the context of the semiregularity map studied there.

The main result of this paper is that the classical Hochschild-Kostant-Rosenberg-Quillen 
theorem on the Hodge decomposition of Hochschild cohomology in the affine case
extends to the global situation, namely that  in the derived category $D(X)$ there 
is a canonical isomorphism of graded $\calox$--algebras 
$$
\bbbh_{X/Y}\cong\bbbs(\bbbl_{X/Y}[1])\,,
$$
where $\bbbl_{X/Y}$ denotes the cotangent complex of $X$ over $Y$. 
In case of smooth morphisms, for instance, $\bbbl_{X/Y}$ is represented by the module of differential forms $\Omega^1_{X/Y}$ and so this returns for (a complex of) $\calox$--modules 
$\calm$ and for any integer $n\in \bbbz$, decompositions \`a la Hodge
\begin{align*}
\HH^n_{X/Y}(\calm)& \cong \bigoplus_{p+q=n}H^{q}(X,\calm\otimes\calt^{p}_{X/Y})\\
\HH_n^{X/Y}(\calm)& \cong 
\bigoplus_{p+q=n}H^{-q}(X,\calm\otimes \Omega^p_{X/Y})\,,
\end{align*}
where $\Omega_{X/Y}^{p}$ or $\calt_{X/Y}^{p}$ denote the usual sheaves of
exterior powers of differential forms or vectorfields, respectively.

For an account of these decompositions in the absolute case of a quasiprojective smooth
variety over a field of characteristic zero, see \cite{Sw} or also \cite{Cal1, Cal2}. 

The interesting new feature in our proof is that the decomposition is achieved through 
characteristic classes, namely (divided and signed) powers of the Atiyah class. 
This allows, for instance, to fit the semiregularity map, as introduced and 
studied in \cite{BFl}, into a larger context.

Let us outline our ideas in the affine case of a ring homomorphism $\Lambda \to A$, 
where $\Lambda$ is a $\bbbq$-algebra. Take a DG resolution $R$ of $\Lambda\to A$, 
that is, a graded free $A$-algebra $R=\bigoplus_{i\le 0}R^{i}$ together with an 
$A$-linear differential, of degree $+1$, such that $H^i(R)=0$ for $i\neq 0$ and $H^0(R)=A$. 
If $B$, which replaces the classical Bar resolution, is then a DG resolution of 
$A$ as a graded free $S:=R\otimes_\Lambda R$-module, then the Hochschild complex 
is represented by the complex of $A$-modules $\bbbh_{B/A}:=B\otimes_S A$, 
and the cotangent complex is represented by the complex 
$\bbbl_{A/\Lambda}:=\Omega^1_{R/\Lambda}\otimes_RA$, see \cite[2.38]{BFl2}.

To deduce the decomposition theorem, note that $B$ is free over $S$, 
whence it is a symmetric algebra over some graded free $S$--module $E$ 
concentrated in negative degrees. Denote $R\otimes 1\subseteq S$ the copy of 
$R$ in $S$ that arises from the inclusions into the left (tensor) factor. 
As $E$ is free over $S$, it admits a connection 
$\nabla:E\to E\otimes_S\Omega^1_{S/R\otimes 1}$ that extends uniquely to an 
$R\otimes 1$-linear derivation 
$\nabla_B:B\to B\otimes_S\Omega^1_{S/R\otimes 1}$ on the algebra $B$. 
The main ingredient in our approach is the {\em Atiyah class\/} 
$\At_B: B\to B\otimes_S\Omega^1_{S/R\otimes 1}[1]$ of this connection, 
given by the $S$--linear commutator $[\partial, \nabla_B]$ with the differential 
$\partial$; see \cite{BFl}. This is an $S$-linear derivation, and we will show that 
Atiyah classes of arbitrary $A$-modules can be obtained from this particular one 
by a suitable tensor product, so we call $\At_B$ the {\em universal\/} Atiyah class; 
see Proposition \ref{4.2.5}.

Consider now the symmetric algebra $\Omega$ over the DG $S$-module 
$\Omega^1_{S/R\otimes 1}[1]$ that is again a DG algebra in a natural way. 
The Atiyah class has a unique $\Omega$-linear extension to a derivation
$\tAt$ on $ B\otimes_S\Omega$ that is locally nilpotent of degree 
zero, whence the map $\exp(-\tAt)$ is an automorphism of DG $\Omega$-algebras. 
Restricting to $B=B\otimes 1$ and tensoring from the right with $R$ 
yields a DG algebra homomorphism $\Phi:B\otimes_SR\to R\otimes_S\Omega$. 
The main step in the proof of the decomposition theorem is to show that this map is 
indeed a quasiisomorphism, see Theorem \ref{main2}, even an isomorphism 
for a suitable choice of $B$, see Theorem \ref{key result}. 
Tensoring further with the augmentation $R\to A$ gives the desired quasiisomorphism 
$$
B\otimes_S A\cong (B\otimes_{S}B)\otimes_{B\otimes 1}A 
\xto{\Phi\otimes 1} (B\otimes_S\Omega)\otimes_{B\otimes 1}A 
\cong A\otimes_S\Omega\,,
$$
since its source is by definition the Hochschild complex $\bbbh_{A/\Lambda}$;
see  \cite[Sect.1.3]{BFl2}; whereas the target is isomorphic to the symmetric 
$A$--algebra over the shifted cotangent complex of $A$ over $\Lambda$, as
$A\otimes_S\Omega^1_{S/R\otimes 1}[1]\cong A\otimes_R\Omega^1_{R/\Lambda}[1]
\cong \bbbl_{A/\Lambda}[1]$. 

To globalize these arguments, we use in the analytic case resolvents on
simplicial schemes of Stein compact sets as developed in \cite{Pa, Fl, BFl, BFl2}. 
Replacing complex spaces by schemes and simplicial schemes of Stein compacta 
by simplicial schemes of affine schemes, one obtains the same results for morphisms 
of schemes over a field of characteristic zero.

We note that in the meantime F.~Schuhmacher  \cite{Sch}, familiar with but independent of an early version of this work, has used the machinery of \cite{BKo} to define Hochschild cohomology for morphisms of complex spaces and algebraic varieties and to give another
proof of the HKR-decomposition theorem.

In the first section we set the stage and define the Atiyah class and the Atiyah-Chern 
character of an augmented analytic algebra. We then formulate and prove the fundamental decomposition theorem in Section 2. Our key applications are to Hochschild theory: in Section 3, we apply the setup of the first section to define the {\em universal\/} Atiyah class and Atiyah-Chern character of a morphism $f:X\to Y$ of complex spaces or schemes over a field of characteristic zero. We show that the universal Atiyah class induces the Atiyah class of any complex of $\calox$--modules. The universal Atiyah-Chern character then yields not only the decomposition of the Hochschild complex and, accordingly, of Hochschild (co-)homology, in view of the decomposition theorem, but we show further that it represents the characteristic homomorphism from Hochschild cohomology to the graded centre of the derived category. In the final Section 4, we fit the {\em semiregularity map\/} for perfect complexes into the picture
and show how it can be be interpreted in terms of the decomposition theorem.

We will use throughout the notions as developed in \cite{BFl2}.

\section{Atiyah Class and Atiyah-Chern Character of an Augmented Algebra}

The key ingredient for our results here is the {\em  Atiyah class\/} of an 
augmented algebra that we introduce first, after some preliminary remarks on 
{\em analytic symmetric algebras\/}.

\subsection{Analytic symmetric algebras}
\begin{sit}\label{d.1}
Let $\calm$ be a coherent $\calox$--module on a complex space $X$, 
and denote $V:=\bbbv(\calm)\to X$ the associated linear fibre space.
If we consider $X$ as a subspace of $V$ via the embedding as the zero section, 
then the symmetric algebra $\bbbs(\calm)$ embeds into the topological restriction 
$\calov|X$ as the subsheaf of functions that are polynomial along the fibres. 
In the following we will call the sheaf
$$
\ts(\calm):=\calov|X
$$
the {\em analytic symmetric algebra\/} of $\calm$. Similarly, if
$K\subseteq X$ is a Stein compact set and $\calm$ a coherent
$\calo_K:=\calox|K$--module, then $\calm$ extends to some neighbourhood
of $K$ as an $\calox$--module, thus we can form $\ts(\calm)$ by
restricting the analytic symmetric algebra of this extended module to
$K$.

For instance, if $\calm$ is a free module on generators $t_1,
\ldots, t_n$, then the stalk $\ts(\calm)_x$, $x\in K$, is just the free
analytic algebra $\calo_{K,x}\{t_1,\ldots, t_n\}$ over
$\calo_{K,x}$ whereas $\bbbs(\calm)_x$ is the subalgebra
$\calo_{K,x}[t_1,\ldots, t_n]$.
\end{sit}

\begin{sit}\label{d.2}
Let $A$ be a simplicial scheme of finite dimension over a set
$I$ and let $W_{*}:=(W_\alpha)_{\alpha\in A}$ be a simplicial system of
Stein compact sets (see \cite[section 2]{BFl}). For a coherent $\calows$--module $\calms$
we can form $\ts(\calms)$ by defining it as before on each simplex.

More generally, assume that $\calrs$ is a DG algebra over $\calows$
and that $\calms$ is a DG $\calrs$--module concentrated in degrees $\le 0$.
The symmetric algebra functor extends to the differential graded context;
see, for example, \cite{Qui2} or \cite{Ill}.
We define then the analytic (differential graded) symmetric algebra to be
$$
\ts(\calms):= \bbbs(\calms)\otimes_\calows \ts(\calm_{*}^{0})\,.
$$
Note that this differs from the usual (differential graded) symmetric algebra only if
$\calms$ contains elements of degree $0$.

Finally, we will call an  $\calrs$--algebra $\calss$ {\em analytic} if
it is isomorphic to a quotient of some analytic symmetric algebra over
$\calrs$ by a homogeneous ideal. Note that then $(W_{*}, \cals_{*}^0)$ may
be considered as a simplicial system on Stein compact subsets.
\end{sit}

\begin{sit}\label{d.3}
Recall that an {\em augmented DG $\calrs$--algebra} is an
$\calrs$--algebra $\calss$ together with a morphism of DG $\calrs$--algebras
$\mu: \calss\to\calrs$, called the {\em augmentation}. 
A morphism of augmented DG $\calrs$--algebras
$\vp:\cals_{*}\to\cals_{*}'$ is a morphism of DG $\calrs$--algebras
that is compatible with the augmentations.

For instance, the analytic symmetric algebra $\ts_{\calrs}(\calms)$
over  a graded DG $\calrs$--module $\calms$ is in a natural way an
augmented DG $\calrs$--algebra; here the augmentation
$\tilde\bbbs_{\calrs}(\calms)\to\calrs$ is the unique homomorphism of
$\calrs$--algebras that 
maps $\calms$ to zero.
\end{sit}

\subsection{The Atiyah class of an augmented algebra}

Let $\calss$ be an augmented DG $\calrs$--algebra that is
analytic as an $\calrs$--algebra. In \cite[section 3]{BFl} we introduced the Atiyah class 
of an $\calss$-module $\calms$ that is bounded above with coherent cohomology as 
an element\footnote{We adopt the general convention that $\Ext^{i}_{\calss}(\calm,\caln) := 
\Hom_{D(\calss)}(\calm,\caln[i])$ for objects $\calm,\caln$ from the indicated derived category.} 
of $\Ext^0_{\calss}(\calms, \calms\dotimes_{\calss}\Omega^1_{\calss/\calrs}[1])$, 
that is, as a morphism 
$$
\At_\calms:\calms\lto \calms\dotimes_{\calss} \Omega^1_{\calss/\calrs}[1]
$$
in the derived category $D(\calss)$. Let us recall in brief how to construct this class.

A DG $\calss$--module $\calms$ that is bounded above with coherent cohomology 
admits a {\em projective approximation\/} $\calps \to \calms$ by \cite[2.19]{BFl}, 
and by  \cite[3.1]{BFl}, the DG module $\calps$ can be endowed with an $\calrs$--linear 
connection 
\begin{align*}
\nabla &: \calp_{*} \to \calp_{*} \otimes_{\calss} \Omega^1_{\calss/\calrs}\,.
\intertext{
With $\partial$ denoting indiscriminately the differential on either source or target,}
[\partial, \nabla] &: \calp_{*} \to \calp_{*} \otimes_{\calss} \Omega^1_{\calss/\calrs}[1]
\end{align*}
is then an $\calss$--linear morphism of DG $\calss$--modules.

The resulting class $\At_{\calms}$ of $[\partial, \nabla]$ in
$$
\Ext^1_{\calss} (\calp_{*} , \calp_{*} \otimes_{\calss}\Omega^1_{\calss /\calrs})
\cong 
\Ext^1_{\calss} (\calm_{*} , \calm_{*} \dotimes_{\calss}\Omega^1_{\calss /\calrs})
$$
is called the {\em Atiyah class\/} of $\calms$; see 3.10 and the remark 
preceding it in \cite{BFl}. 

We need the following simple observation; see \cite[3.3]{BFl}. It uses the canonical 
isomorphisms $\bbbs^n(\Omega^1_{\calss/\calrs}[1])\cong \Omega^{n}_{\calss/\calrs}[n]$; 
see \cite{Ill}; and yields that the (powers of the) Atiyah class of such a DG module are 
well defined in the derived category.

\begin{lem}\label{lem BFl}
Let $f:\calm'_{*}\to\calms$ be a morphism of DG $\calss$--modules. If
$$
\nabla': \calm'_{*} \to \calm'_{*}\otimes_{\calss} \Omega^1_{\calss /\calrs}
\quad\text{and}\quad
\nabla: \calms \to \calms\otimes_{\calss} \Omega^1_{\calss /\calrs}
$$
are connections, then $(f\otimes 1)\circ[\partial,\nabla']^{n}$ and
$[\partial,\nabla]^{n}\circ f$ are {\em homotopic\/} for every integer $n\ge 0$; 
they represent the same class in
$H^{n}(\Hom_{\calrs} (\calm'_{*}, \calm_{*} \otimes_{\calss}\Omega^{n}_{\calss/\calrs}))$. 

Thus, if $\nabla: \calms \to \calms\otimes_{\calss}\Omega^1_{\calss/\calrs}$ is any
connection and one of the modules $\calms$, $\Omega^1_{\calss/\calrs}$ is {\em locally free\/} 
over $\calss$, then the powers of the deduced Atiyah class of $\calms$ are well 
defined and are represented in the derived category by the morphisms 
$$
[\partial, \nabla]^{n}: \calms \to 
\calms\otimes_{\calss} \bbbs^n(\Omega^1_{\calss/\calrs}[1])=
\calms\dotimes_{\calss}\bbbs^n(\Omega^1_{\calss/\calrs}[1])\,.
$$
\qed
\end{lem}

Applying this to $\calms=\calrs$, viewed as an $\calss$-module via the augmentation, 
we make the following definition. 

\begin{defn}
\label{d.5}
The Atiyah class of $\calrs$ as an $\calss$-module,
$$
\Atu:=\At_\calrs:\calrs\lto \calrs\dotimes_{\calss} \Omega^1_{\calss/\calrs}[1]
$$ 
in $D(\calss)$, will be called the {\em universal Atiyah class} of the augmented algebra. 
\end{defn}

This notation will be justified in Proposition \ref{4.2.3}. It will be important to have 
explicit descriptions of this Atiyah class and we next provide two of those.

\begin{sit}\label{d.4}
 There exists
a free DG algebra resolution $\nu:\calbs\to \calrs$ of $\calrs$ as an
$\calss$--algebra; see, for example, \cite[1.1.2(1)]{BFl2}. We may 
and will assume $\calbs$ to be generated in degrees $<0$, whence
$\calbs\cong \bbbs_{\calss}(\cales)$ for some graded  free 
DG $\calss$--module $\cales$. By the results of
\cite[Sect. 3.1]{BFl}, the $\calss$--module $\cales$ admits an $\calrs$--linear 
connection $\nabla:\cales\to\cales\otimes_{\calss}\Omega^1_{\calss/\calrs}$. 
Such a connection extends uniquely to an $\calrs$--linear
{\em derivation\/} on the $\calss$--algebra $\calbs$,
$$
\nabla_{\calbs}: \calbs\to \calbs\otimes_{\calss}\Omega^1_{\calss/\calrs}\,,
$$
such that $\nabla_{\calbs}$ restricts to the exterior differential 
$d:\calss\to \Omega^1_{\calss/\calrs}$ on $\calss$ and to the given connection 
$\nabla$ on  $\cales\subseteq \calbs$. 
Note that $\nabla_{\calbs}$ is then in particular an  $\calrs$--linear connection. 
By definition, $\Atu$ is then represented by the commutator
$$
\At_\calbs:= [\partial, \nabla_{\calbs}]:\calbs\to
\calbs\otimes_{\calss}\Omega^1_{\calss/\calrs}[1].
$$
Observe that, by construction, this morphism is even an $\calss$--linear derivation.   
\end{sit}

\begin{sit}\label{d.6}
Using the Zariski-Jacobi sequence, one can give yet another description of
this Atiyah class. The DG algebra homomorphisms $\calrs\to 
\calss \to \calbs$ determine the Zariski-Jacobi sequence
\begin{equation}
     \label{eq:zj}
     0\to  \calbs  \otimes_{\calss}\Omega^{1}_{\calss/\calrs}
     \xto{d\vp} \Omega^{1}_{\calbs/\calrs} \to
     \Omega^{1}_{\calbs/\calss}\to 0
\end{equation}
that is an exact sequence of DG $\calbs$--modules as $\calbs$ 
is free over $\calss$, thus, $\Omega^{1}_{\calbs/\calss}$ projective in the
category of DG $\calbs$--modules. 
Completing this exact sequence to a distinguished 
triangle in the derived category of DG $\calbs$--modules gives 
rise to a morphism
$$
\alpha:\Omega^{1}_{\calbs/\calss}\lto
(\calbs \otimes_{\calss}\Omega^{1}_{\calss/\calrs})[1] \,.
$$
Exploiting once again that  $\Omega^{1}_{\calbs/\calss}$ is projective,
one may choose $\alpha$ as an actual morphism of complexes that is then 
unique up to homotopy. Composing $\alpha$ with the universal $\calss$--derivation
$$
d: \calbs\lto \Omega^{1}_{\calbs/\calss}
$$
results in an $\calss$--linear derivation of DG $\calbs$--modules
$$
\alpha \circ d: \calbs\lto
(\calbs \otimes_{\calss}\Omega^{1}_{\calss/\calrs})[1] \,.
$$
We claim that $-\alpha \circ d$ represents the Atiyah class of $\calss/\calrs$. In
fact, as the module $ \Omega^{1}_{\calbs/\calss}$ 
in the sequence (\ref{eq:zj}) is projective, we can
find a morphism of graded $\calbs$--modules 
$\beta:\Omega^{1}_{\calbs/\calrs}\to  \calbs
\otimes_{\calss}\Omega^{1}_{\calss/\calrs}$ splitting that
sequence. The map $-\alpha$ is then the map induced by $[\partial,
\beta]$. Moreover, the map
$$
\nabla=\beta \circ d:\calbs\lto  \calbs
\otimes_{\calss}\Omega^{1}_{\calss/\calrs}
$$
is an $\calrs$--linear  connection on $\calbs$ as an $\calss$--module. The 
claim follows then from the equalities
$
-\alpha \circ d=[\partial,\beta]d= [\partial, \beta \circ d]= [\partial,\nabla]
$, with the second one due to the fact that $d$ is a morphism of DG modules.
\end{sit}

\subsection{The Atiyah-Chern character of an augmented algebra}
\begin{sit}
\label{Atu}
In the following, it will be convenient to extend the Atiyah class 
$\At_{\calbs}$ from \ref{d.4} to a {\em derivation\/} on the larger DG algebra
$\calbs\otimes_{\calss}\Omegas$, where we set
$$
\Omegas:=\bbbs_{\calss}(\Omega^1_{\calss/\calrs}[1])\,.
$$
More precisely, $\At_\calbs$ extends uniquely to an $\Omegas$-linear derivation 
$$
\tAt: \calbs\otimes_{\calss}\Omegas\to \calbs\otimes_{\calss}\Omegas
\quad\mbox{via}\quad
\tAt(b\otimes\omega):=\At_\calbs(b)\cdot (1\otimes \omega)
$$ 
for local sections $b$ in $\calbs$ and $\omega$ in $\Omegas$.  
The following lemma exhibits the pertinent properties of this extension.
\end{sit}

\begin{lem}
\label{lem main2} 
With notation as just introduced, the following hold.
\begin{enumerate}[\quad\rm (1)]
\item 
\label{main2.1}
The morphism $\tAt$ is an $\Omegas$-linear derivation of degree $0$.

\item
\label{main2.2} 
The derivation $\tAt$ is locally nilpotent\footnote{This means that for
every local section, say, $h$ of  $\calbs \otimes_{\calss}\Omegas $ there is an integer
$n=n(h)\in \bbbn$ such that $\tAt^n(h)=0$.}.

\item
\label{main2.3} 
The exponential of the derivation $-\tAt$ induces an automorphism of DG
$\Omegas$--algebras
$$
\exp(-\tAt):\calbs\otimes_{\calss} \Omegas \lto
\calbs\otimes_{\calss} \Omegas\,.
$$
\end{enumerate}
\end{lem}

\begin{proof}
(\ref{main2.1}) holds by construction.
To deduce (\ref{main2.2}), write $\cales=\bigoplus_{i\le -1}\cale_{*}^{(i)}$, where
$\cale_{*}^{(i)}$ is a graded free module over $\calss$ with
generators only in degree $i$. Let $\calb_{*}^{(k)}$
denote the $\calss$--subalgebra of $\calbs$ generated by 
$\bigoplus_{i>k}\cale_{*}^{(i)}$. 
As $\nabla$ is homogeneous of degree zero, it restricts to a 
connection on $\bigoplus_{i>k}\cale_{*}^{(i)}$,
whence $\nabla_\calbs(\calb_{*}^{(k)})\subseteq \calb_{*}^{(k)}$, and so $\tAt$
restricts for each $k$ to a derivation
$$
\tAt:  \calb_{*}^{(k)}\otimes_{\calss}\Omegas
\lto  \calb_{*}^{(k)}\otimes_{\calss}\Omegas\,.
$$
We show by decreasing induction on $k$ that this restriction is
locally nilpotent. As $\tAt$ vanishes on $1\otimes\Omegas$, this is
true for $k=0$. Assume that the claim is shown for the level
$k$. As $\calb_{*}^{(k-1)}\otimes_{\calss}\Omegas$ is generated 
as a $\calb_{*}^{(k)} \otimes_{\cals}\Omegas$--algebra by
$\cale_{*}^{(k)}$, it is sufficient to verify that a given
local section $e$ of  $\cale_{*}^{(k)}$ is annihilated by a sufficiently
high power of $\tAt$. Now the elements $\partial(e)$ as well as 
$\partial\nabla(e)$ are contained in the subalgebra 
$\calb_{*}^{(k)}\otimes_{\calss}\Omegas$
whence $\tAt( e)=\At_\calbs(e) \in \calb_{*}^{(k)}\otimes_{\calss}\Omegas$. Applying
the induction hypothesis, (\ref{main2.2}) follows. 
The final assertion (\ref{main2.3}) is a well known consequence of 
(\ref{main2.2}); see, for example, \cite{Ren}.
\end{proof}

\begin{sit}
\label{d6.55}
Restricting the automorphism $\exp(-\tAt)$, as just constructed in 
\ref{lem main2}(\ref{main2.3}), to the
subalgebra $\calbs\cong \calbs\otimes 1$ of
$\calbs\otimes_{\calss}\Omegas$, we obtain a homomorphism of 
DG $\calss$--algebras
$$
\AC_\calbs:=\exp(-\tAt)|_{\calbs}:\calbs \lto \calbs\otimes_{\calss} \Omegas\,,
$$
explicitly, for a local section $b$ of $\calbs$,
\begin{align*}
\AC_\calbs(b)&= \exp(-\tAt(b\otimes 1)) = 
\sum_{n\ge 0}\frac{(-1)^{n}}{n!} \At_{\calbs}^{n}(b)\,,
\end{align*}
where the powers of 
$\At_{\calbs}(b)\in \calbs\otimes\Omega^{1}_{\calss/\calrs}$ occurring in the last term 
are calculated in $\calbs\otimes_{\calss} \Omegas$. In this sense, $\AC_\calbs=\exp(-\At_\calbs)$. 
\end{sit}

In light of this, we make the following definition.

\begin{defn}
\label{d6.6} 
Let $\calss$ be an augmented analytic $\calrs$--algebra and 
$\nu:\calbs\to\calrs$ a DG algebra resolution of $\calrs$ over $\calss$. The
{\em universal Atiyah-Chern character\/} 
of $\calss/\calrs$ is the element
$$
\AC_u:=\AC_\calbs\in 
\Ext^0_{\calss}(\calbs ,\calbs\otimes_{\calss}\Omegas)
\cong
\Ext^0_{\calss}(\calrs ,\calrs\dotimes_{\calss}\Omegas)\,.
$$
\end{defn}

The class of $\AC_\calbs$ in $D(\calss)$ is independent of the choice of a resolution $\calbs$ 
and we can summarize the construction and properties established so far as follows.

\begin{prop} \label{uniqueness}
The universal Atiyah-Chern character of an augmented analytic algebra
$\calss/\calrs$ is a uniquely defined morphism 
$AC_u:\calrs\to \calrs\dotimes_{\calss}\Omegas$ in the derived 
category of DG $\calss$--modules.
\qed
\end{prop}

\section{The Decomposition Theorem for Augmented Algebras}
We now turn to our first main result that says that the Atiyah-Chern character 
of an augmented analytic algebra $\calss/\calrs$ induces an algebra {\em isomorphism\/} 
in the derived category of DG $\calrs$--modules
$$
\Phi_\calrs:
\calrs\dotimes_{\calss}\calrs \lto \calrs\dotimes_{\calss}\Omegas
$$
provided that $\calss$ is a free analytic $\calrs$--algebra.

\subsection{Statement of the key results}
Assume that $\calss$ is a free analytic $\calrs$-algebra so that 
$\Omegas =\bbbs_{\calss}(\Omega^{1}_{\calss/\calrs}[1])$ is, in particular, 
a free $\calss$-module and the Atiyah-Chern character $\AC_{u}$ in $D(\calss)$ 
is represented by a morphism of DG $\calss$--modules 
$\calrs\to \calrs\otimes_{\calss}\Omegas$ .
Tensoring it with $\calrs$ in $D(\calrs)$ gives a morphism
$$
\AC_u\dotimes_{\calss}\id: \calrs \dotimes_{\calss}\calrs
\lto (\calrs\otimes_{\calss}\Omegas)\dotimes_{\calss}\calrs\,.
$$
Composing this with the natural morphism 
$$
(\calrs\otimes_{\calss}\Omegas)\dotimes_{\calss}\calrs\lto 
(\calrs\otimes_{\calss}\Omegas)\otimes_{\calss}\calrs   \cong \calrs\otimes_{\calss}\Omegas
$$
returns in $D(\calrs)$ a morphism
$$
\Phi_\calrs:  \calrs \dotimes_{\calss}\calrs
\lto \calrs\otimes_{\calss}\Omegas\,.
$$
\begin{sit}
\label{stage}
To represent as well $\Phi_\calrs$ by an actual morphism of DG modules, we return to the 
homomorphism of  DG $\calss$--algebras
$$
\AC_\calbs:\calbs \lto \calbs\otimes_{\calss} \Omegas
$$
from \ref{d6.55}. The composition of DG $\calss$--module homomorphisms
\bdi
\calbs &\rTo^{\AC_\calbs} &  \calbs\otimes_{\calss} \Omegas
&\rTo^{\nu\otimes \id} & 
\calrs\otimes_{\calss} \Omegas\,,
\edi
has as its target a DG $\calrs$--module, thus, factors uniquely through a morphism of
DG $\calrs$--modules
$$
\Phi_\calbs:  \calbs \otimes_{\calss}\calrs
\lto \calrs\otimes_{\calss}\Omegas\,.
$$
As $\AC_u$ is represented by $\AC_\calbs$, it follows that
$\Phi_\calbs$ 
represents $\Phi_\calrs$ in $D(\calrs)$.
\end{sit}

The key result is now the following decomposition theorem.

\begin{theorem}\label{main2}
Let $\calss/\calrs$ be an augmented analytic algebra 
such that $\calss$ is a free analytic algebra over $\calrs$.
The homomorphism of DG $\calrs$--algebras $\Phi_\calbs$ from above is then a
quasiisomorphism and so the morphism $\Phi_\calrs$ is an isomorphism in $D(\calrs)$.
\end{theorem}

We will in fact show the following stronger statement.

\begin{theorem}
     \label{key result}
Under the assumptions of \ref{main2}, one can construct a graded free 
DG $\calss$--algebra resolution $\nu: \calbs\to\calrs$ such that the morphism
$$
\Phi_\calbs: \calbs\otimes_{\calss}\calrs \lto
\calrs\otimes_{\calss} \Omegas\,.
$$
is an {\em isomorphism} of DG $\calrs$--algebras.
\end{theorem}

Note that this result implies \ref{main2} immediately.
The remainder of this section is devoted to a proof of \ref{key result}, 
by constructing explicitly the algebra resolution with the required property.
To this end, we first recall the construction and structure
of acyclic DG algebras in characteristic zero, as the explicit contracting 
homotopy of such algebras in terms of the inverse to the Euler derivation 
takes centre stage in the proof.

\subsection{Acyclic algebras}
To begin, we recall the construction in characteristic zero 
of the acyclic $\calrs$--algebra defined by a DG $\calrs$--module.  
This tool is due to Quillen \cite[App.II]{Qui2} in case of DG algebras 
over a field; its classical precursor is the explicit contraction of the 
Koszul complex on the variables over a polynomial algebra; see, 
for example, \cite[X.206, Exerc.~\S 9.2]{Alg}.  

\begin{sit}
\label{acyclic algebra}
Let $\calgs$ be a DG $\calrs$--module generated in degrees $\le
0$ and set $T\calgs:= \calgs[1]$, a shifted copy of $\calgs$ 
with the inherited differential. 
We equip the DG module $T\calgs\oplus\calgs$ with the 
differential
$$
\partial (Tg+g') := -T\partial (g) +g+\partial (g')
\leqno(*)
$$
for local sections $g,\, g'$ of $\calgs$. With this
differential, $T\calgs\oplus\calgs$ is an acyclic DG
$\calrs$--module; indeed, it is the mapping cone over the identity
map on $\calgs$. An explicit contracting homotopy is given by the
$\calrs$-linear map $\delta:T\calgs\oplus\calgs\to T\calgs\oplus\calgs$  
of degree $-1$ defined through $\delta(Tg+g'):= Tg'$.

In view of \ref{d.1}, the (graded and analytic) symmetric algebra
$\calas:=\tilde\bbbs_{\calrs}(T\calgs\oplus\calgs)$ on
$T\calgs\oplus\calgs$ is in a natural way an augmented DG
$\calrs$--algebra. Its augmentation ideal is denoted
$$
\calj_{*}:=\tilde
\bbbs^{>0}_{\calrs}(T\calgs\oplus\calgs)=
(T\calgs\oplus\calgs)\calas\subseteq \calas\,.
$$
The differential $\partial$ on $\calas$ is uniquely characterized
by the property that it restricts to the given differential on $\calrs$ 
and that $\partial |(T\calgs\oplus\calgs)$ is the differential
described in $(*)$. Moreover, the homotopy $\delta$ extends uniquely to an
$\calrs$--linear derivation of degree $-1$ that we denote again
$\delta:\calas\to\calas[-1]$.
\end{sit}

\begin{defn}
\label{def:1}
The augmented DG algebra $(\calas, \partial)$ as just
constructed will be called the {\em acyclic DG $\calrs$--algebra
on\/} $\calgs$.
\end{defn}

To justify the terminology, consider the commutator
$$
\epsilon := [\partial,\delta]:\calas\lto\calas\,.
$$
The following lemma identifies $\epsilon$ as the {\em Euler
derivation\/} with respect to the ``symmetric degree'' on
$\calas$.

\begin{lem}
\label{euler}
The derivation $\epsilon$ on $\calas$ is $\calrs$--linear
of degree $0$.  It has the following properties.
\begin{enumerate}[\quad\rm (1)]
\item
     \label{mult}
The restriction of $\epsilon$ to $\bbbs^d(T\calgs\oplus\calgs)$ 
is multiplication by $d$.  In particular, $\epsilon$ maps the 
augmentation ideal $\calj_{*}$ of $\calas$ bijectively to itself.
\item
     \label{commutes}
The derivation $\epsilon$ commutes with $\partial$ and
with $\delta$.  If $\xi:\calj_{*}\to\calj_{*}$ denotes the inverse
of $\epsilon$ on $\calj_{*}$, then $\xi$ commutes with
$\partial$ and $\delta$ as well.
\end{enumerate}
\end{lem}

\begin{proof}
The restriction of $\epsilon$ to $\calrs$ is equal to the
commutator of $\partial|\calrs$ and $\delta|\calrs=0$ and so
is identically zero on $\calrs$. Thus $\epsilon$ is an
$\calrs$-linear derivation of degree 0  with
$\epsilon|(T\calgs\oplus\calgs)=\id$, whence (1) follows
from
$\calrs$--linearity and the product rule.

In view of $ \partial^2=\delta^2=0$, the first part of
(2) is a consequence of the graded Jacobi identity, whereas the second
part follows from the first.
\end{proof}

\begin{cor}
\label{contraction}
Structure morphism and augmentation of the augmented DG
$\calrs$--algebra $(\calas,\partial)$ are homotopy
equivalences; equivalently, the augmentation ideal 
$\calj_{*}=\ker(\calas\to\calrs)$
is contracted by the homotopy $h:=\xi\delta=\delta\xi$.
\end{cor}

\begin{proof}
Indeed, the map $\xi\delta$ of degree $-1$ is a contracting 
homotopy: on $\calj_{*}$, one has
$$
[\partial, \xi\delta] = \xi[\partial,\delta]
= \xi\epsilon = \id_{\calj_{*}}
$$
by \ref{euler}(\ref{commutes}).
\end{proof}

\subsection{The key construction}
\begin{sit}
     \label{special}
We apply the foregoing construction in the following special case.
Consider a presentation of the given augmented graded free
$\calrs$--algebra $\calss$, say
$$
\calss\cong\tilde\bbbs_{\calrs}(\calfs)\,,
$$
where $\calfs=\bigoplus_{i\le 0}\calf_{*}^{(i)}$ is a graded free
$\calrs$--module with the direct summand $\calf_{*}^{(i)}$ freely
generated in degree $i$. 
Recall that $\mu:\calss\to\calrs$ denotes the given augmentation, 
and let $j:\calrs\to\calss$ be the structure homomorphism of the 
$\calrs$--algebra $\calss$. 

Replacing the $\calrs$-submodule $\calfs\subseteq\calss$ by the 
isomorphic copy $(\id_{\calss}-j\circ \mu)(\calfs)\subseteq\calss$, 
we may assume that the augmentation ideal is generated by
$\calfs\subseteq\calss$, so that
$$
\cali_{*}:=\ker (\mu) =
\calfs\calss\subseteq\calss
$$
and the given differential $\partial$ on $\calss$ maps
$\calfs$ into the augmentation ideal $\cali_{*}$.  

We now trivialize the differential on $\calfs$ as follows.  Let
$\tilde\calfs$ be a second copy of $\calfs$ as graded
$\calrs$--module, where the identification
$\calfs\cong\tilde\calfs$ is written $f\mapsto \tilde f$.  Endow
$\tilde\calfs$ with the unique DG $\calrs$--module structure
under which all local module generators are cycles and let
$(\calbs,\partial)$ be the acyclic $\calrs$--algebra on
this new DG $\calrs$--module $\tilde\calfs$ as in \ref{def:1}.  
We still denote $\delta$ the canonical $\calrs$--linear derivation
of degree $-1$ on $\calbs$, let $\epsilon = [\partial,\delta]$
be the corresponding Euler derivation and $\xi$ its inverse on the
augmentation ideal $\calj_{*}$ of $\calbs$. 
Recall that then $h=\xi\delta =\delta\xi$ contracts $\calj_{*}$
by \ref{contraction}.
\end{sit}

The key step in proving theorem \ref{key result} is now the next
proposition, where we endow $\calbs$ with the structure of an
$\calss$-algebra, that is, embed $\calss$ as an augmented 
DG subalgebra of $\calbs$. The idea is straightforward:
local generators of degree 0 are necessarily cycles, they can be mapped
to the corresponding generators in $\tilde\calfs$. 
We then proceed by decreasing induction on degrees, using that we already
will know where to map the images under the differential in $\calss$ 
and then applying the homotopy $h$ on $\calbs$ to find a suitable image 
in $\calbs$ for the generators of lower degree. 
The gist of the proof is to control the inductive process. 

\begin{prop}
     \label{key}
With notations as introduced above, there is a morphism of
augmented DG
$\calrs$--algebras $\vp: \calss \to
\calbs$ that is uniquely determined through the
recursion
$$
\vp(f) = \tilde f +h\circ\varphi\circ\partial(f)
$$
for every local $\calrs$--generator $f$ of $\calfs$.  It satisfies
the following properties.
\begin{enumerate}[\quad\rm (1)]
\item
The morphism $\vp$ realizes $\calbs$ 
as the graded free $\calss$--algebra
$$
\calbs\cong\bbbs_{\calss}(T \tilde\calfs\otimes_\calrs\calss)\,.
$$ 

\item
     \label{factors}
The image of $\epsilon \circ \varphi:\calss\to \calbs$
is contained in $\cali_{*}\calbs:= \vp(\cali_{*})\calbs$.
\end{enumerate}
\end{prop}

\begin{proof}
Let 
$$
\cals_*^{(k)}:=\bbbs_{\calrs}\left(\bigoplus_{i >
k}\calf_{*}^{(i)}\right)\subseteq\calss\quad \mbox{and}\quad 
\calb_*^{(k)}:=\bbbs_{\calrs}\left(\bigoplus_{i >
k}(T{\tilde \calf}_{*}^{(i)}\oplus {\tilde
\calf}_{*}^{(i)})\right)\subseteq\calbs
$$
denote the DG $\calrs$--subalgebras of $\calss$ and
$\calbs$, respectively, that are generated by elements of degree greater 
than $k$, so that $\calss = \liminj_{k}\cals_{*}^{(k)}$ and 
$\calbs = \liminj_{k}\calb_{*}^{(k)}$. These algebras are again augmented
$\calrs$--algebras, and clearly $\calb_*^{(k)}$ is the acyclic
$\calrs$--algebra on $\bigoplus_{i > k}\tilde\calf_{*}^{(i)}$. In
particular, the morphisms $\epsilon,\delta$ stabilize 
$\calb_{*}^{(k)}$, and $\xi,h$ are defined on its augmentation
ideal $\calj_{*}$.  We denote by
$\cali_{*}^{(k)} :=\cali_{*}\cap \cals_{*}^{(k)} 
=\ker(\cals_{*}^{(k)}\to\calrs)$ 
the augmentation ideal of $\cals_{*}^{(k)}$. 

To establish existence and claimed properties of $\varphi$, it
suffices to construct, by descending induction on $k$, 
morphisms of augmented DG $\calrs$--algebras 
$\varphi^{(k)}:\cals_{*}^{(k)}\to
\calb^{(k)}_{*}$ that satisfy the following properties.
\begin{enumerate}
\item[(1)$_k$]\quad 
$\varphi^{(k)}|\cals_{*}^{(k+1)}=\varphi^{(k+1)}$,
and $\varphi^{(k)}$ realizes $\calb^{(k)}_{*}$ as
$$
\calb^{(k)}_{*}\cong
\left(\calb^{(k+1)}_{*}\otimes_{\cals_{*}^{(k+1)}}
\cals_{*}^{(k)}\right)\otimes_{\calrs}
\bbbs_{\calrs}(T\tilde\calf^{(k+1)}_{*})\,.
$$

\item[(2)$_k$]\quad 
$\epsilon \circ \varphi^{(k)}:
\cals^{(k)}_{*}\to \calb_{*}$ maps $\cals^{(k)}_{*}$ into
$\cali^{(k)}_{*}\calb_{*}:=\vp^{(k)}(\cali^{(k)}_{*})\calb_{*}\subseteq
\calj_{*}$.
\end{enumerate}
Once the family $(\varphi^{(k)})_{k}$ has been constructed, $\vp =
\liminj_{k}\varphi^{(k)}$ will satisfy claim (1) in view of (1)$_k$,
whereas claim (2) will follow from the conditions (2)$_k$.

Before we turn to the construction of the morphisms $\varphi^{(k)}$,
let us prove the following claim that we will need in the course of
the proof:
\begin{itemize}
     \item[($\dagger$)$_{k}$] If the preceding two conditions
(1)$_k$, (2)$_k$
     are satisfied for $\varphi^{(k)}$, then $\cali^{(k)}_{*}\calb_{*}$ is
     mapped bijectively onto itself under $\epsilon$.
\end{itemize}
Indeed, $\epsilon$ is a derivation, hence using the product rule and
$(2)_k$ it follows that $\epsilon$ maps
$\cali^{(k)}_{*}\calb_{*}$ to itself.  By \ref{euler}, $\epsilon$ is
bijective and semisimple on $\calj_{*}$ as a $\bbbc$--linear
endomorphism, and this implies that its restriction to any vector
subspace stable under $\epsilon$ is bijective.  As $\varphi^{(k)}$
respects the augmentation, we have
$\cali^{(k)}_{*}\calb_{*}\subseteq \calj_{*}$, and as
$\cali_{*}^{(k)}\calbs$ was just seen to be stable under $\epsilon$,
the claim ($\dagger$)$_{k}$ follows.

Now we construct the morphisms $\varphi^{(k)}$. Clearly, for $k\ge
0$, we set $\cals_{*}^{(k)} = \calrs =\calb_{*}^{(k)}$ and
$\cali_{*}^{(k)} = 0$, and let $\varphi^{(k)}$ be the identity
map, whence the claimed properties hold vacuously true.  
Let us assume now that for some $k\le 0$ a morphism 
$\varphi^{(k)}$ has already been constructed
that satisfies the required conditions.

For a local section $f$ in $\calf^{(k)}_{*}$ of degree $k$ note that
$\partial f$ is a local section in $\cals^{(k)}_{*}$, so that
$\varphi^{(k)}\partial (f)$ is already defined as a local section of
$\calb^{(k)}_{*}$.  Indeed, $\partial f$ is even a local section in
$\cali^{(k)}_{*}$, thus, $\varphi^{(k)}\partial (f)$ is in
$\cali^{(k)}_{*}\calb_{*}\subseteq \calj_{*}$, and so $\xi$
is defined on it.  Accordingly,
\begin{align}
\label{eq:1star}
\tag{$*$}
\varphi^{(k-1)}(f):= \tilde f + h\varphi^{(k)}\partial(f) = 
\tilde f +\delta\xi\varphi^{(k)}\partial(f)
\end{align}
is a well defined local section in $\calb^{(k-1)}_{*}$.  As
$\cals^{(k-1)}_{*}$ is the symmetric $\cals_{*}^{(k)}$--algebra on
$\calf^{(k)}_{*}\otimes_{\calrs}\cals^{(k)}_{*}$, this map extends
uniquely to a morphism of $\cals^{(k)}_{*}$--algebras $\varphi^{(k-1)}:
\cals^{(k-1)}_{*} \to\calb^{(k-1)}_{*}$.  It follows directly from this
explicit description that (1)$_{k-1}$ holds. As $\delta$, $\partial$ and
$\xi$ stabilize $\calj_{*}$ and as $\tilde f\in \calj_{*}$ for $ f\in
\cali_{*}$, it is immediately seen from $(*)$ that 
$\varphi^{(k-1)}$ is compatible with the augmentation.  

Now we verify that $\varphi^{(k-1)}$ is a morphism of DG algebras.  
First observe that for a local section $f$ as before
$$
\partial\xi\varphi^{(k)}\partial(f)=
\xi\partial\varphi^{(k)}\partial(f)=
\xi\varphi^{(k)}\partial^{2}(f)=0\,,
$$
where we have used \ref{euler}(\ref{commutes}) and the fact that
$\varphi^{(k)}$ is a DG morphism. In particular, this shows that
$\delta\partial\xi\varphi^{(k)}\partial(f) =0$, and thus
\begin{equation}
\label{eq:2star}
\tag{$**$}
\partial\delta\xi\varphi^{(k)}\partial(f) =
[\partial,\delta]\xi\varphi^{(k)}\partial(f) =
\varphi^{(k)}\partial(f)\,,
\end{equation}
as $\epsilon = [\partial,\delta]$ and $\xi$ are inverse to each other
on $\varphi^{(k)}\partial (f)\in \cali^{(k)}_{*}\calbs\subseteq
\calj_{*}$.  
Next apply $\partial$ to equation (\ref{eq:1star}) to obtain
\begin{alignat*}{2}
\partial\varphi^{(k-1)}(f)&=\partial\delta\xi\varphi^{(k)}\partial
(f)&\quad\quad&\text{as $\partial {\tilde f}=0$\,,}\\
&= \varphi^{(k)}\partial(f)&&\text{by (\ref{eq:2star})\,,}\\
&= \varphi^{(k-1)}\partial(f)&&\text{by (1)$_{k-1}$\,,}
\end{alignat*}
whence $\varphi^{(k-1)}$ is indeed a morphism of DG
$\calrs$--algebras.

It remains to verify (2)$_{k-1}$.  In view of (2)$_{k}$ and
(1)$_{k-1}$ it is enough to show that $\epsilon\varphi^{(k-1)}(f)$ is
in $\cali^{(k-1)}_{*}\calbs$ for each local section $f$ as before.  To
this end, note first that $\varphi^{(k-1)}(f)$ is in
$\cali^{(k-1)}_{*}\calb_{*}$ by definition, whence it suffices to show the current  
claim for $\epsilon\varphi^{(k-1)}(f)-\varphi^{(k-1)}(f)$.  Applying
$(\epsilon-\id)$ to $(*)$ and using that $\epsilon(\tilde f)=\tilde f$;
see  \ref{euler}(\ref{mult}); we obtain
\begin{alignat*}{2}
(\epsilon-\id)\varphi^{(k-1)}(f)&=
(\epsilon-\id)\delta\xi\varphi^{(k)}\partial (f)
&\qquad&\\
&=\delta\xi(\epsilon-\id)\varphi^{(k)}\partial (f)
&&\text{by \ref{euler}(\ref{commutes}),}\\
&=\delta(\id - \xi)\varphi^{(k)}\partial(f)
&&\text{as $\xi\epsilon = \id$ on $\calj_{*}$.}
\end{alignat*}
We will now show that this term is indeed already in 
$\cali_{*}^{(k)}\calbs\subseteq \cali_{*}^{(k-1)}\calbs$.
Even more generally, we verify that $\delta(\id - \xi)\varphi^{(k)}(g)$
is contained in $\cali^{(k)}_{*}\calb_{*}$ for any local section
$g\in\cals^{(k)}_{*}$. Note that the term in question is 
certainly in $\calj_{*}\cap\calb_{*}^{(k)}$, as $\varphi^{(k)}(g)$ is in
$\calb_{*}^{(k)}$ and $\delta\circ(\id-\xi)$ maps $\calb_{*}^{(k)}$ to
$\calj_{*}\cap\calb_{*}^{(k)}$. As $\epsilon$ is bijective on 
$\cali_{*}^{(k)}\calb_{*}$ in view of ($\dagger$)$_k$, it suffices now 
to show that
\begin{align}
\label{eq:3stars}
\tag{$*{*}*$}
\epsilon\delta(\id - \xi)\varphi^{(k)}(g) =
\delta\epsilon(\id - \xi)\varphi^{(k)}(g) =
\delta(\epsilon - \id)\varphi^{(k)}(g)
\end{align}
is contained in $\cali^{(k)}_{*}\calbs$.  Suppose that this is
satisfied for homogeneous local sections $g_{1}, g_{2}$ of
$\cals^{(k)}_{*}$.  We claim that it is then as well satisfied for
$g=g_{1} g_{2}$.  In fact, expanding by means of the product rule for
the derivations $\delta$ and $\epsilon$ yields
\begin{align*}
\delta(\epsilon - \id)\varphi^{(k)}(g_{1} g_{2})&=
\left(\delta(\epsilon - \id) \varphi^{(k)}(g_{1})\right)
\varphi^{(k)}(g_{2})\\
&\hphantom{=} +(-1)^{|g_{1}|}
\varphi^{(k)}(g_{1})\left(\delta(\epsilon - \id)
\varphi^{(k)}(g_{2})\right)\\
&\hphantom{=} +(-1)^{|g_{1}|}\epsilon\varphi^{(k)}(g_{1})
\delta\varphi^{(k)}(g_{2})\\
&\hphantom{=} + \delta\varphi^{(k)}(g_{1})
\epsilon\varphi^{(k)}(g_{2})\,.
\end{align*}
The first two terms on the right-hand side are in $\cali^{(k)}_{*}\calb_{*}$ by assumption,
whereas the remaining two are so by (2)$_{k}$.  This leaves us to establish
that the element from (\ref{eq:3stars}) is in $\cali^{(k)}_{*}\calb_{*}$
for a local generator $g$ of $\calfs$ of degree greater than $k$.
In this situation, however, $\varphi^{(k)}(g) = \tilde g
+\delta\xi\varphi^{(k)}\partial(g)$ by (1)$_{k}$ and (\ref{eq:1star}),
whence the element in (\ref{eq:3stars}) is equal to
$$
\delta\epsilon(\tilde g+ \delta\xi \varphi^{(k)}\partial(g))
-\delta(\tilde g+\delta\xi \varphi^{(k)}\partial(g)) = 0\,,
$$
as $\epsilon(\tilde g)=\tilde g$, the derivation $\epsilon$
commutes with $\delta$, and $\delta^2=0$. The proof of the proposition
is thus complete.
\end{proof}

Next we describe the action of an extended Atiyah class 
on the resolution constructed in the foregoing proposition, and we keep
the notation from there.

\begin{con}
\label{convention}
To get a simple explicit description of the Atiyah class $\At_{\calbs}$
as introduced in \ref{d.5}, it is convenient to suppress in the following 
the morphism $\vp:\calss\to \calb_{*}$, that is,
$\calss$ will be considered as a DG $\calrs$--subalgebra of
$\calbs$. Accordingly, a local section $f$ of
$\calfs\subseteq\calss$ becomes identified with the element
$\tilde f+h\circ \partial(f)$ of $\calbs$.
\end{con}

\begin{sit}
\label{Atu2}
By construction, 
$\calbs=\bbbs_\calrs(T\tilde\calfs\oplus \tilde\calfs)\cong\bbbs_{\calss}(\cales)$ 
with $\cales:=T\tilde\calfs\otimes_\calrs\calss$. 
With notation as in \ref{Atu},  let us take the connection $\nabla$ on the 
$\calss$-module $\cales$ with $\nabla|T\tilde\calfs=0$. It extends as in 
\ref{Atu} to a unique derivation $\nabla_\calbs:\calbs\to\calbs\otimes_{\calss}\Omegas$
that restricts to the exterior differential $d$ on $\calss$. The next Lemma describes 
how the corresponding Atiyah class $\At_\calbs=[\partial,\nabla_\calbs]$ acts 
on the algebra generators of $\calbs$.

\begin{lem}
\label{3.2.10} 
With $\nabla_\calbs$ as above, the Atiyah class
$$
\At_\calbs=[\partial,\nabla_\calbs]:
\calbs\lto\calbs\otimes_{\calss}\Omegas
$$
is the $\calss$-linear derivation uniquely characterized 
by\footnote{Note that $\Omega^1_*=\Omega^1_{\calss/\calrs}[1]$, so to keep track of the shift we will write $Tdf$ instead of $df$.}
$$
\At_\calbs(T\tilde f )=-1\otimes Tdf +\nabla_\calbs h\partial(f)
$$
for a local section $f$ of $\calfs$.
\end{lem}

\bproof
First recall from \ref{acyclic algebra}($*$) that $\partial(T\tilde f)
= -T(\partial\tilde f) + \tilde f$. Now $\partial(\tilde f)=0$ in $\calbs$
by definition of that algebra, and $\tilde f= f - h\partial(f)$ by the 
convention \ref{convention} just made above. As $\nabla_\calbs(T\tilde f)=0$
by choice of the connection, we get in summary
$\At_\calbs(T\tilde f ) = -\nabla_\calbs\partial(T\tilde f) = -\nabla_\calbs(f - h\partial(f))$,
which expands to the claimed form as $\nabla_\calbs(f) =1\otimes Tdf$.
\eproof

As in the proof of \ref{key}, write $\calf=\bigoplus_{i\le
0}\calf^{(i)}$, with $\calf^{(i)}$ a free $\calrs$--module 
generated in degree $i$. Different from there, let now
$\calb_{*}^{(k)}\subseteq\calbs$  be the $\calss$--subalgebra
of $\calbs$ generated by $\bigoplus_{i>k}T\tilde \calf_{*}^{(i)}$,
and introduce $\Omega_{*}^{(k)}\subseteq
\Omegas$ as the $\calss$--subalgebra generated by
$\bigoplus_{i>k}Td(\calf_{*}^{(i)})$, where $d$ denotes the exterior
differential on $\Omegas$ as in \ref{Atu}. 
In view of Lemma \ref{3.2.10} above, $\At_\calbs$ restricts to a derivation $$
\At_\calbs:\calb_{*}^{(k)}\lto \calb_{*}^{(k)}\otimes\Omegas\,.
$$
\end{sit}

\subsection{Proof of the decomposition theorem}
Now we are ready to prove the decomposition theorem in its strong form.

\bproof[Proof of Theorem \ref{key result}]
As noted above, $AC_\calbs=\exp(-\At_\calbs)$ maps 
$\calb_{*}^{(k)}$ into
$\calb_{*}^{(k)}\otimes_{\calss}\Omega_*^{(k)}\subseteq
\calb_{*}\otimes_{\calss}\Omega_*^{(k)}$. Combining this with the construction in \ref{stage}, we see that $\Phi_\calbs$, as defined there,
restricts for each integer $k$ to a morphism
\begin{align}
\Phi_\calbs^{(k)}:\tilde\cala_{*}^{(k)}:=
\calb_{*}^{(k)}\otimes_{\calss}\calrs
\lto \cala_{*}^{(k)}:=
\calrs\otimes_{\calss}\Omega_{*}^{(k)} 
\end{align}
of DG $\calrs$--algebras, and, in turn,
$\Phi_\calbs= \liminj_{k} \Phi_\calbs^{(k)}$.
It thus suffices to verify, by decreasing induction on $k$, that each
$\Phi_\calbs^{(k)}$ is an isomorphism.

For $k=0$ this is evident as $\calb_{*}^{(0)}=\Omega_*^{(0)}=\calss$
and the morphism reduces to the identity.
Now assume that $\Phi_\calbs^{(k)}$ is already an isomorphism for some $k\le 0$. 
The algebra $\tilde\cala_{*}^{(k-1)}$ is a symmetric $\tilde\cala_{*}^{(k)}$--algebra
over  the free module $T\tilde \calf_{*}^{(k)}\otimes_{\calrs} \tilde\cala_{*}^{(k)}$ and,
similarly, $\cala_{*}^{(k-1)}$  is a symmetric
$\cala_{*}^{(k)}$--algebra over the free module
$Td(\tilde\calf_{*}^{(k)})\otimes_{\calrs}\cala_{*}^{(k)}$, and the hypothesis
is that $\tilde\cala_{*}^{(k)}\cong \cala_{*}^{(k)}$ under $\Phi_\calbs^{(k)}$.

If $f$ is a local $\calrs$-generator of
$\calf_{*}^{(k)}$, then by \ref{3.2.10}, the Atiyah class
satisfies $\At_\calbs( T\tilde f) =
-1\otimes Tdf +\beta$ for some $\beta\in \cala_{*}^{(k)}$. As $\Phi_\calbs^{(k)}$ is induced by 
$\exp(-\At_\calbs)$ on $\calb_{*}^{(k)}\otimes_{\calss}1$, this gives
$$
\Phi_\calbs^{(k)}( T\tilde f\otimes 1)=(\mbox{residue class of } \exp(-\At_\calbs)( T\tilde f) \mbox{ in }\calrs\otimes_{\calss}\Omega_*)=
1\otimes Tdf +\gamma
$$
for some  $\gamma\in \cala^{(k)}$.  We are thus in the situation as expressed
through the commutative diagram
\bdi
{\bbbs_{ \tilde\cala_{*}^{(k)}}}(T\tilde \calf_{*}^{(k)}{\otimes_{\calrs}} \tilde\cala_{*}^{(k)})
&=\tilde\cala_{*}^{(k-1)}&\rTo^{\Phi_\calbs^{(k-1)}}& \cala_{*}^{(k-1)}&
={\bbbs_{ \cala_{*}^{(k)}}}(Td(\tilde\calf_{*}^{(k)}){\otimes_{\calrs}}\cala_{*}^{(k)})\\
&\uTo&&\uTo\\
&\tilde\cala_{*}^{(k)}&\rTo^{\Phi_\calbs^{(k)}}_{\cong}&\cala_{*}^{(k)}
\edi
and the morphism on top is a sum $\vp_{1}+\vp_{0}$, where 
\begin{align*}
\vp_{1}: T\tilde \calf_{*}^{(k)}{\otimes_{\calrs}} \tilde\cala_{*}^{(k)}\xto{\ \cong\ } 
Td(\tilde\calf_{*}^{(k)}){\otimes_{\calrs}}\cala_{*}^{(k)}
\end{align*}
sends $T\tilde f\otimes 1$ to $1\otimes Tdf$, and is thus an 
$\tilde\cala_{*}^{(k)}$--isomorphism of DG modules over the isomorphism 
$\Phi_\calbs^{(k)}$ of DG algebras, while $\vp_{0}$, that sends $T\tilde f\otimes 1$
to $\gamma\in \cala_{*}^{(k)}$, is some homomorphism over $\Phi_\calbs^{(k)}$.
The proof is now completed by applying the simple lemma hereafter 
to each stalk of this diagram.
\eproof

\begin{lem}
Let $\tilde A\to A$ be an isomorphism of graded rings, 
$\tilde F$ a graded $\tilde A$--module, and $F$ a graded $A$--module. 
Let $\vp:\bbbs_{\tilde A}(\tilde F)\to \bbbs_A(F)$ be a morphism of graded
$\tilde A$--algebras such that $\vp|\tilde F=\vp_0+\vp_1$, where
$\vp_0:\tilde F\to A$ and
$\vp_1:\tilde F \to F$ are linear over the given ring isomorphism. 
If $\vp_1$ is an isomorphism of modules, then $\vp$ is an isomorphism.
\end{lem}

\begin{proof}
Clearly we may assume that $\tilde A=A$. The morphism
$\psi=-\vp_0\circ \vp_1^{-1}+\vp_1^{-1}:F\to\bbbs_A(\tilde F)$ 
induces then a map $\bbbs_A(F)\to\bbbs_{\tilde A}(F)$ that is an 
inverse to $\vp$.
\end{proof}

\begin{rem}
The quasiisomorphism of DG $\calrs$--modules $\calbs\otimes_{\calss}\calrs$
to $\bbbs_\calrs(\calrs\otimes_{\calss}\Omega^1_{\calss/\calrs})$ 
induced by the Atiyah-Chern character is well defined up to those homotopies 
that arise from different choices of the connection on $\calbs$. 
In view of \ref{Atu2}, such different choices of connections
reflect simply different choices of the presentation 
$\calss\cong\bbbs_{\calrs}(\calfs)$ of $\calss$ as free $\calrs$--algebra.

We note now that one can use \ref{key} more directly to construct a 
quasiisomorphism of DG $\calrs$--algebras in the opposite direction, 
$\bbbs_\calrs(\calrs\otimes_{\calss}\Omega^1_{\calss/\calrs})
\xto{\simeq} \calbs\otimes_{\calss}\calrs$. 
However we are not able to confirm that this
morphism is well defined in the derived category, in particular, 
we fail to know whether it represents an inverse quasiisomorphism
to the isomorphism induced by the Atiyah-Chern character.

The construction of this quasiisomorphism is as follows, 
where we use freely the notation from \ref{key}. 
First note that $ \calbs\otimes_{\calss}\calrs\cong\calbs/\cali_{*}\calbs$,
where $\cali_{*}$ is the augmentation ideal of $\calss$.
Composing the $\calrs$--linear derivation $\delta\varphi:\calss\to
\calb_{*}$ of degree $-1$ with the natural $\calss$--algebra
morphism $\calb_{*}\to \calb_{*}/\cali_{*}\calb_{*}$ induces an
$\calrs$--linear derivation, say,
$\overline{\delta\vp}:\calss\to \calb_{*}/\cali_{*}\calb_{*}$ of
degree $-1$.  Now \ref{key}(\ref{factors}) yields that the commutator
$$
[\partial, \delta\varphi] = \partial(\delta\varphi) +
(\delta\varphi)\partial =
\partial\delta\varphi + \delta\partial\varphi = \epsilon\circ\varphi
$$
maps $\calss$ to $\cali_{*}\calb_{*}$, whence
$\overline{\delta\vp}$ is compatible with the differentials.
Hence $\overline{\delta\vp}$ is a DG derivation, thus
represented by an $\calss$--linear map of DG modules
$h:\Omega^1_{\calss/\calrs}[1]\to \calb_{*}/\cali_{*}\calb_{*}$
of degree 0 that further factors through an
$\calrs$--linear map
$\bar{h}:\calrs\otimes_{\calss}\Omega^1_{\calss/\calrs}[1]
  \to \calb_{*}/\cali_{*}\calb_{*}$. The
algebra map induced by $\bar h$,
$$
\bbbs_{\calrs}(\bar{h}): \bbbs_{\calrs}(\calrs\otimes_{\calss}
\Omega^1_{\calss/\calrs}[1]) \lto\calb_{*}/\cali_{*}\calb_{*}
$$
is then a morphism of augmented DG $\calrs$--algebras. Let us show
that this map is an isomorphism. As
$\calss\cong\bbbs_{\calrs}(\calfs)$, the universal
derivation restricts to an isomorphism of
$\calrs$--modules $\bar d=d\otimes 1: \calfs[1]\xto{\cong}
\Omega^1_{\calss/\calrs}[1] \otimes_{\calss}\calrs$, and
for a local generator $f$ in $\calfs$ the element
$\bar{h}\bar df$ is by construction the residue class of
$T\tilde f$ in $\calb_{*}/\cali_{*}\calb_{*}$.  Thus,
$\calb_{*}/\cali_{*}\calb_{*}$ being a free $\calrs$--algebra over
$T\tilde\calfs\cong \bar{h}\bar d(\calfs[1])$, the
algebra morphism induced by $\bar{h}$ is an isomorphism of augmented
DG $\calrs$--algebras, as required.
\end{rem}

\section{Decomposition of the Hochschild Complex}

Let $X\to Y$ be a fixed morphism of complex spaces\footnote{We 
leave again to the reader the simpler case of morphisms of schemes 
over a field of characteristic zero}. In this part we will
apply the decomposition theorem \ref{main2} to obtain a natural
quasiisomorphism from the Hochschild complex
$\bbbh_{X/Y}$ into the symmetric algebra $\bbbs(\bbbl_{X/Y}[1])$ over
the shifted cotangent complex.

\subsection{The decomposition of the Hochschild complex for 
morphisms of complex spaces}
\begin{sit}\label{4.1.1}
Let $\fX=(X_{*},W_{*},\calrs,\calbs)$ be a free extended resolvent for
$X$ over $Y$, with $\calss$ the enveloping algebra of $\calrs$ and 
$\nu:\calbs\to\calrs$ a free DG algebra resolution of $\calrs$ as $\calss$--algebra
as in \cite{BFl2}. This means in particular that $\calbs$ is the symmetric 
algebra over a graded free $\calss$--module, say, $\cales$ over $\calss$. 
Restricting $\calss$ topologically to $W_{*}$ we obtain an analytic free 
DG $\calrs$--algebra $\calss|W_{*}$ that can be written as the analytic 
symmetric algebra of some free module $\calfs$ over $\calrs$, see \ref{d.1}. 
We view $\calss$ always as an augmented $\calrs$--algebra via the structure 
morphism $j= j_{1}:\calrs\to \calrs\otimes 1\subseteq \calss$ and the augmentation 
given by the mutiplication morphism $\mu:\calss\to \calrs$. 
\end{sit}

\begin{sit}
\label{universal.1}
We are now able to apply the decomposition theorem \ref{main2}. 
The hypotheses of that result are satisfied in this context by definition of an 
extended free resolvent, and so we obtain an $\calrs$-linear quasiisomorphism
$$
\Phi_\calbs: \calbs\otimes_{\calss}\calrs \lto
\calrs\otimes_{\calss} \bbbs_{\calss}(\Omega^1_{\calss/\calrs}[1])\cong 
\bbbs_\calrs(\Omega^{1}_{\calrs/Y}[1])\,.
$$
Using the quasiisomorphism $\calrs\to\caloxs$ that is part of the data of an 
extended resultant, we obtain a quasiisomorphism of DG $\caloxs$--modules
\begin{align}
\label{at*}
\Phi_{X_*}:=
\Phi_\calbs\otimes_{\calrs}\id_\caloxs: \bbbh_{X_{*}/Y}=
\calbs\otimes_{\calss}\caloxs \lto
\bbbs_\caloxs(\Omega^1_{\calrs/Y}[1]\otimes_\calrs\caloxs)\,.
\end{align}
By definition, $\Omega^1_{\calrs/Y}\otimes_\calrs\caloxs$ represents the 
cotangent complex $\bbbl_{X_{*}/Y}$ of $X_{*}$ over $Y$, whence the algebra 
on the right is the symmetric algebra over the shifted cotangent complex 
$\bbbl_{X_{*}/Y}[1]$. 
Taking finally the associated \v{C}ech complexes yields the following 
decomposition theorem for the Hochschild complex of a morphism of 
complex spaces.
\end{sit}

\begin{theorem}\label{main4}
Let $X\to Y$ be a morphism of complex spaces. In the derived
category $D(X)$, the universal Atiyah-Chern character yields a canonical 
isomorphism of graded $\calox$--algebras 
$$
\Phi=\Phi_{X/Y}:= \calcb(\Phi_{X_*}): 
\bbbh_{X/Y}\lto\bbbs(\bbbl_{X/Y}[1])\,,
$$
where $\bbbl_{X/Y}$ denotes the cotangent complex of $X$ over $Y$ and
$\bbbs(\bbbl_{X/Y}[1])$ the derived symmetric algebra; cf.\
\cite[4.2.8(3)]{BFl2} or \cite{Ill}. In particular, there are natural decompositions
of Hochschild (co-)homology in terms of the derived exterior powers of the 
cotangent complex,
\begin{eqnarray*}
\HH^n_{X/Y}(\calm)&\cong &
\bigoplus_{p+q=n}\Ext^q_X(\Lambda^p(\bbbl_{X/Y}), \calm)\\
\HH_n^{X/Y}(\calm)& \cong &
\bigoplus_{p+q =n}\Tor_q^X(\Lambda^p(\bbbl_{X/Y}), \calm)\,.
\end{eqnarray*}
If $\calm=\calox$ then the second\footnote{but not the first!} of these isomorphisms is compatible
with the algebra structures.
\qed
\end{theorem}

In case of smooth morphisms, for instance, this yields decompositions \`a la Hodge,
in the following sense.

\bcor
For a smooth morphism of complex spaces $X\to Y$ of relative dimension $d$, 
the universal Atiyah-Chern character effects natural decompositions
\begin{align*}
\HH^n_{X/Y}(\calm)&\cong \Ext^n_{X\times_Y X}(\calox,\calm)\cong
\bigoplus_{p+q=n}H^{q}(X,\calm\otimes\calt^{p}_{X/Y})\\
\HH_n^{X/Y}(\calm)&\cong \Tor^{X\times_Y X}_n(\calox,\calm) \cong 
\bigoplus_{p+q=n}\Tor_q^X(\Omega^p_{X/Y}, \calm)\,,
\end{align*}
where $\Omega_{X/Y}^{p}$ or $\calt_{X/Y}^{p}$ denote the usual $\calox$--modules of
exterior powers of differential forms or vectorfields, respectively. 

If $\calm=\calox$, then the second of these isomorphisms is compatible
with the algebra structures.

In particular, for $\calm=\calox$ the structure sheaf, respectively for
$\calm=\omega_{X/Y}\cong \Omega^{d}_{X/Y}$ the relative canonical sheaf,
one obtains
\begin{align*}
\HH_n^{X/Y}(\calox)&\cong \bigoplus_{p+q=n}H^{-q}(X,\Omega^p_{X/Y})\\
\HH^n_{X/Y}(\omega_{X/Y})&\cong\bigoplus_{p+q=n}H^{q}(X,\Omega^{d-p}_{X/Y}) \,.
\end{align*}
\qed
\ecor
As mentioned in the introduction, a detailed account of these 
decompositions in the absolute case of a quasiprojective smooth
variety over a field of characteristic zero, can be found in \cite{Sw} 
or also \cite{Cal1,Cal2}.

\begin{rem}\label{4.1.3}
Using the functoriality of Atiyah classes it follows that the
decomposition of the Hochschild complex is compatible with
morphisms. More precisely, given a commutative diagram of complex
spaces
\bdi[small]
X'&\rTo & X\\
\dTo &&\dTo\\
Y'&\rTo& Y
\edi
the corresponding diagram
\bdi[small]
Lf^*(\bbbh_{X/Y})&\rTo^\cong & Lf^*(\bbbs(\bbbl_{X/Y}[1]))\\
\dTo && \dTo\\
\bbbh_{X'/Y'}&\rTo^\cong & \bbbs(\bbbl_{X'/Y'}[1])
\edi
commutes as well, where the horizontal maps are given by the
exponentials of the Atiyah classes and the vertical arrows are the
canonical ones. This follows easily from the fact that the Atiyah
class and its powers are compatible with morphisms, see
e.g.\ \cite[3.14]{BFl} for a related situation.
\end{rem}

\subsection{Applications to augmented complex spaces}
Let $X$ be a complex space and $Y$ a closed subspace. Assume 
that there is a holomorphic map $f:X\to Y$ restricting to the identity 
on $Y$. Using the decomposition theorem we can show the following 
result.

\bpro
\label{augmented}
The universal Atiyah-Chern character of the subspace $Y$ in $X$ 
provides a ca\-non\-ical isomorphism 
$$
\caloy\dotimes_\calox\caloy\cong \bbbs_Y(\bbbl_{Y/X})\,.
$$
In particular, for every $\caloy$-module $\calq$ and each integer $p$,
we have isomorphisms
$$
\ba{rcl}
\Tor_{p}^X(\caloy, \calq)&\cong &
\Tor_{p}^Y(\bbbs_Y(\bbbl_{Y/X}), \calq)\\
\Ext^{p}_X(\caloy, \calq)&\cong &
\Ext^{p}_Y(\bbbs_Y(\bbbl_{Y/X}), \calq)\,.
\ea
$$ 
\epro

Before proving this result we note the following important special case. 

\bcor
If $Y$ is locally a complete intersection in $X$ with normal bundle $\caln_{Y/X}$, then 
$$
\caloy\dotimes_\calox\caloy\cong \bbbs_Y(\caln_{Y/X}^\vee)\,.
$$
In particular, for every $\caloy$-module $\calq$ and each integer $p$, 
we have isomorphisms
$$
\ba{rcl}
\Tor_{p}^X(\caloy, \calq)&\cong &
H^{-p}(Y,\bbbs_Y(\caln_{Y/X}^\vee)\otimes_\caloy \calq)\\
\Ext^{p}_X(\caloy, \calq)&\cong &
H^{p}(Y,(\bbbs_Y(\caln_{Y/X})\otimes_\caloy \calq)\,.
\ea
$$ 
\qed
\ecor

\bproof[Proof of Proposition \ref{augmented}]
Let $(X_{*},W_{*},\calrs)$ be a resolvent of $X\to Y$. Intersecting the simplicial 
space $X_{*}$ with $Y$ gives a simplicial space $Y_{*}$ that is the nerve 
of a covering of $Y$. We let $\calozs$ be the topological 
restriction of $\caloxs$ to $Y_{*}$. The map $X\to Y$ induces a morphism 
$\caloys\to \calozs$, and the embedding $Y\hto X$ a surjection 
$\calozs\to \caloys$, so that $\caloys\to \calozs\to \caloys$ is the identity. 
The topological restriction, say, $\calss$ of $\calrs$ to $Y_{*}$ is a free 
$\caloys$-algebra that admits an augmentation 
$\calss\to\caloys$ given by the composition 
$\calss\to\calozs\to\caloys$. Using the decomposition theorem there is 
thus a natural decomposition 
$$
\caloys\dotimes_{\calss}\caloys
\cong \bbbs(\Omega^1_{\calss/Y}[1]\otimes _{\calss}\caloys)
\cong \bbbs(\bbbl_{X_{*}/Y_{*}}[1]\dotimes_\caloxs\caloys)\,.
$$
As $\calss$ is quasiisomorphic to $\calozs=\caloxs|Y_{*}$, the first term is 
just the topological restriction of  $\caloys\dotimes_\caloxs\caloys$ to $Y_{*}$. 
Taking  \v{C}ech complexes yields an isomorphism 
$$
\caloy\dotimes_\calox\caloy\cong \bbbs(\bbbl_{X/Y}[1]\dotimes _\calox\caloy)\,.
$$
Finally, using the distinguished triangle 
$$
\bbbl_{X/Y}\dotimes _\calox\caloy \to \bbbl_{Y/Y}\to \bbbl_{Y/X} 
\to \bbbl_{X/Y}[1|\dotimes _\calox\caloy
$$
associated to the maps $Y\to X\to Y$, and noting that $\bbbl_{Y/Y}=0$,
we find $\bbbl_{Y/X} \cong \bbbl_{X/Y}[1]\dotimes _\calox\caloy$, whence the 
result follows.
\eproof

\begin{rem} 
In the absolute affine case, where $Y=\Spec K$ is a field of characteristic 
zero and $X=\Spec A$ a noetherian augmented $K$--algebra, \ref{augmented} yields the
decomposition $K\dotimes_{A}K\cong \bbbs_{K}(\bbbl_{K/A})$, a result due to Quillen
\cite{Qui3}. In this case, however, the result has a generalization to augmented algebras 
over {\em any field\/}, except that the symmetric algebra has to be replaced by a divided 
power algebra and $\bbbl_{K/A}$ by the dual of the {\em homotopy Lie algebra\/} of $A$ 
over $K$; see \cite{Av} for an account of this. The decomposition theorem thus extends 
certainly beyond characteristic zero, but the precise category of morphisms of (analytic) 
spaces  to which it generalizes and which proper form such generalization should take 
are currently not known.
\end{rem}

\section{The Universal Atiyah Class and Atiyah-Chern Character}
In this section we show that the universal Atiyah class of a morphism $f:X\to Y$
of complex spaces induces the Atiyah classes of all complexes in the derived category 
of $X$, whence the name. As a byproduct, we extend the results from \cite{BFl} to show
that the Atiyah classes and the Atiyah-Chern character are well defined for {\em any\/} 
complex, not just for those locally bounded above with coherent cohomology.
We first describe the situation on a fixed extended resolvent and then descend to $X$ itself
using again the \v{C}ech construction.

\subsection{Atiyah classes of $\calrs$--modules versus the universal 
Atiyah class} 
Let $D^{-}_{c}(\calrs)$ denote the derived category of bounded DG 
$\calrs$--modules with coherent cohomology as considered in \cite{BFl}. 
There, we developed the theory of (powers of) Atiyah classes of such modules, see 
also section 1.2.
Here we first verify that those classes are induced by the universal Atiyah class,
before we then extend the results to the full derived category $D(\calrs)$.

\begin{sit}
Consider a free extended resolvent $(X_{*},W_{*}, \calrs,\calbs)$ 
for a morphism of complex spaces $X\to Y$ and a DG $\calrs$--module 
$\calms$ that is bounded above with coherent cohomology. 
To relate $\At_{\calms}$ to the universal Atiyah class
$\Atu$, let us form the enveloping algebra
$$
\calss := \calrs \tilde\otimes_{\calo_Y} \calrs
$$
as in \ref{4.1.1}. This algebra admits two different $\calrs$--algebra 
structures; a ``left'' one and a ``right'' one, given by the subalgebras
$\calrs\otimes 1$, respectively $1\otimes \calrs$ of $\calss$.  The
multiplication map $\mu :\calss\to \calrs$ endows
$\calss$ with the structure of an augmented $\calrs$--algebra, with
respect to either structure, and generally we view $\calss$ as an 
$\calrs$--algebra with respect to the left structure. That said,
we can find, as in section 1.2, a ``universal" connection
$$
\nabla_\calbs:
\calbs\lto\calbs\otimes_{\calss}\Omega^1_{\calss/\calrs\otimes 1}
$$
so that $\nabla_\calbs$ is an $\calrs\otimes 1$--linear
derivation restricting to the exterior differential on $\calss$.
The DG $\calss$ module $\Omega^1_{\calss/\calrs\otimes 1}$ is 
canonically isomorphic to $\calrs\tilde\otimes_{\caloy} \Omega^1_{\calrs/Y}$. 
Employing this isomorphism, $\nabla_\calbs$ can thus be viewed as a derivation
$$
\nabla_\calbs:
\calb_{*}\lto\calb_{*}\otimes_{\calss}
(\calrs\tilde \otimes_{\caloy} \Omega^1_{\calrs/Y})
\cong \calb_{*}\otimes_{1\otimes \calrs}\Omega^1_{\calrs/Y}
$$
that is $\calrs\otimes 1$--linear and restricts to the exterior differential on 
$1\otimes \calrs$.
If now $\calms$ is any DG $\calrs$--module, we can form the tensor product over
$\calrs=\calrs\otimes 1$ to obtain a connection
$$
\id_{\calms}\otimes \nabla_\calbs:
\calms\otimes_{\calrs\otimes 1}\calbs\lto
 \calms\otimes_{\calrs\otimes 1}\calbs\otimes_{1\otimes
\calrs}\Omega^1_{\calrs/Y}
$$
on the $1\otimes\calrs$--module $\calms\otimes_{\calrs\otimes
1}\calbs$. Finally, since $\calbs$ and $\Omega^1_{\calrs/Y}$ are free 
$\calrs$-modules that are bounded above and the given morphism 
$\nu:\calbs\to \calrs$ is a quasiisomorphism, the natural morphisms
\begin{align*}
\id_{\calms}\otimes  \nu&:
\calms\otimes_{\calrs\otimes 1}\calbs
\lto
\calms\otimes_\calrs\calrs \cong \calms\\
\id_{\calms}\otimes \nu\otimes \id &: \calms
\otimes_{\calrs\otimes 1}\calbs
\otimes_{1\otimes\calrs} \Omega^1_{\calrs/Y}
\lto
\calms\otimes_{\calrs}\Omega^1_{\calrs/Y}
\end{align*}
are quasiisomorphism as well; see \cite[2.17 (1)]{BFl}.

Forming $[\partial , \id_{\calms}\otimes \nabla_\calbs]$ yields in view of \ref{lem BFl} 
the following result that justifies to call $\Atu$, or its representative $\At_\calbs$, 
the {\em universal\/} Atiyah class as we proposed in Definition \ref{d.5}.
\end{sit}

\begin{prop}\label{4.2.3}
If $\calms$ is a DG $\calrs$--module in $D^-_c(\calrs)$, then 
$\id_{\calms}\otimes \At_\calbs$, when viewed as a morphism of $1\otimes \calrs$-modules, 
represents $\At_{\calms}$ in $D(\calrs)$. 
Equivalently\footnote{The reader should note that there is no need to derive the 
tensor product as $\Omega^1_{\calrs/Y}$ is {\em flat\/}, even free, as $\calrs$--module. 
Similar remarks explain why we may use ``actual'' tensor products throughout this section.}, 
$$
\At_\calms=\id_{\calms}\otimes_{ \calrs\otimes 1} \Atu:\calms\lto  
\calms\otimes_\calrs\Omega^1_{\calrs/Y}
$$
as morphisms in $D(\calrs)$.
\hfill $\square$
\end{prop}

The reader should note the following subtle points in this proposition. 
The universal Atiyah class $\Atu$ is a morphism of $\calss$-modules, so 
$\id_\calms\otimes\Atu$ is as well a morphism of $\calss$-modules, 
thus, we can view it as a morphism of $\calrs$-modules from the right or of 
$\calrs$-modules from the left, and the resulting morphisms in $D(\calrs)$ are, in general, 
completely different! 
As a morphism of $\calrs\otimes 1$-modules, the map $\id_\calms\otimes \Atu$ is  even 
zero, since  $\id_\calms\otimes 1:\calms \to \calms\otimes_{\calrs\otimes 1}\calbs$ is 
$(\calrs\otimes 1)$-linear, whence $\id_\calms\otimes \Atu$ is represented by the real world 
morphism $m\mapsto m\otimes \At_\calbs(1)=0$ on $\calms$ in $D(\calrs\otimes 1)$. 

Looking at it the other way around, we may use this result to obtain Atiyah
classes for {\em arbitrary} DG modules in $D(\calrs)$. 
More generally, we will use the universal Atiyah-Chern character 
$\AC_{u}$ from \ref{d6.6} to define the Atiyah-Chern character
for any DG $\calrs$--module.

\begin{defn}\label{4.2.4}
For a DG module $\calms\in D(\calrs)$, the morphism
$$
\At_{\calms}\in \Hom_{D(\calrs)}(\calms,
\calms\otimes_{\calrs} \Omega^1_{\calrs/Y}[1])
$$
in $D(\calrs)$ defined through the commutative diagram
\bdi
\calms\otimes_{\calrs\otimes 1}\calb_{*}
&\rTo^{1_{\calms}\otimes\At_\calbs} &
  \calms\otimes_{\calrs\otimes 1}\calb_{*}\otimes_{1\otimes
\calrs}\Omega^1_{\calrs/Y}[1]\\
\dTo^{\id_{\calms}\otimes\nu}_{\simeq} && 
\dTo_{\id_{\calms}\otimes \nu\otimes \id}^{\simeq}\\
\calms & \rTo^{\At_{\calms}:=\id_\calms\otimes \Atu} &
\calms\otimes_{\calrs} \Omega^1_{\calrs/Y}[1]
\edi
will be called the {\em Atiyah class\/} of $\calms$. 

In a similar way, the universal Atiyah-Chern character $\AC_u$ 
as represented by $\At_\calbs$ gives rise to the {\em Atiyah-Chern character\/} of 
$\calms$ the morphism
$$
\AC_{\calms}:=\id_\calms\otimes\AC_u: \calms \lto
\calms\otimes_{\calrs} \Omegas
$$
in $D(\calrs)$, with $\Omegas:=\bbbs(\Omega^1_{\calrs/Y}[1]).$
\end{defn}

\begin{prop}
\label{4.2.5} 
Let $\calm$ be a DG module in $D(\calrs)$.
\begin{enumerate}[\quad\rm (1)]
\item
\label{4.2.5.1}
The Atiyah class $\At_{\calms}$ and the Atiyah-Chern character
$\AC_{\calms}$ are functorial in $\calms$ with respect to morphisms
in $D(\calrs)$.

\item
\label{4.2.5.2} 
If $\calms$ admits a connection
$\nabla:\calms\to \calms\otimes_{\calrs}\Omega^1_{\calrs/Y}$
then the morphism of DG modules $[\partial, \nabla]$ represents 
$$\At_{\calms}\in 
\Ext^0_\calrs(\calms, \calms\otimes_{\calrs}\Omega^1_{\calrs/Y}[1])\,.
$$

\item
\label{4.2.5.3} 
There exists always a quasiisomorphism of DG modules 
$\vp:\calm'_{*}\xto{\simeq} \calms$ such that
$\calm'_{*}$ admits a connection $\nabla'$. The morphism
$(\vp\otimes\id)\circ[\partial,\nabla']\circ \vp^{-1}$ represents then $\At_{\calms}$
in $D(\calrs)$.

\item
\label{4.2.5.4} The Atiyah classes of $\calms$ and of the shifted DG module
$\calms[1]$ are related through $\At_{\calms[1]} = (-\At_{\calms})[1]$.
\end{enumerate}
\end{prop}

\begin{proof}
(\ref{4.2.5.1}) is immediate from the definition, and (\ref{4.2.5.2}) follows from
\ref{lem BFl}  applied to the quasiisomorphism
$\calms\to  \calms\otimes_{\calrs\otimes 1}\calb_{*}$ and the
connections $\nabla$ and $\id_{\calms}\otimes \nablau$. For (\ref{4.2.5.3}),
take $\calm_{*}':=\calms\otimes_{\calrs\otimes 1}\calbs$ and
$\nabla:=\id_{\calms}\otimes \nablau$. The final point (\ref{4.2.5.4})
follows as in \cite[3.7]{BFl}.
\end{proof}

\begin{sit}
To summarize the results of the preceding proposition in an abstract fashion,
disembodied of the DG module $\calms$, let us call a morphism
$\vp:F\to G$ of exact functors between triangulated categories 
{\em exact\/} if, for any distinguished triangle $\Delta$ in the source 
category of $F,G$, the morphism $\vp(\Delta):F(\Delta)\to G(\Delta)$ is a 
morphism of distinguished triangles in the target category. In detail, 
this means the following: as $F,G$ are exact functors, they come equipped 
with isomorphisms of functors $\epsilon:F\circ T\xto{\cong} T\circ F$ and 
$\eta:G\circ T\xto{\cong} T\circ G$, where $T$ denotes the respective 
translation functor. The requirement that $\vp(\Delta)$ is a morphism of triangles 
then simply means that the diagram of morphisms of functors
\bdi
F\circ T&\rTo^{\vp\circ T}&G\circ T\\
\dTo^{\epsilon}&&\dTo_{\eta}\\
T\circ F&\rTo^{T\circ \vp}&T\circ G
\edi
commutes. For the Atiyah class, considered as a morphism $\At_{?}:\id_{D(\calrs)}\to
\id\otimes_{\calrs}\Omega^1_{\calrs/Y}[1]$, this amounts precisely to the
property \ref{4.2.5}(\ref{4.2.5.4}).
\end{sit}

Now we can rephrase \ref{4.2.5} as follows.

\begin{cor}
\label{At transformation}
The Atiyah class defines an exact morphism
\begin{align*}
\At_{?}:\id_{D(\calrs})\lto \id_{D(\calrs)}\dotimes_{\calrs}\Omega^1_{\calrs/Y}[1]\,,
\end{align*}
while the Atiyah-Chern character defines an exact morphism 
\begin{align*}
\AC_{?}:\id_{D(\calrs})\lto \id_{D(\calrs)}\dotimes_{\calrs}\Omegas
\end{align*}
of endofunctors on $D(\calrs)$.\qed
\end{cor}

We finish this subsection with another interpretation of the Atiyah class using
the exact sequence of free DG $\calrs$--modules
$$ 
0 \to \Omega^1_{\calrs/ Y}\cong \cali_{*} / \cali^2_{*} \to \calss /
\cali^2_{*} \to
\calrs \to 0, 
$$
where $\cali_{*} = \Ker \mu$ is the kernel of $\mu$. To fix signs we
identify $\Omega^1_{\calrs/ Y}$ and $ \cali_{*} / \cali^2_{*}$ as in
\cite[3.17]{BFl} via
$da := 1\otimes a-a\otimes 1$. Tensoring the sequence from the left
with $\calms$ over $\calrs=\calrs\otimes 1$ returns an exact sequence of DG 
$\calrs=1\otimes\calrs$--modules
\begin{equation}
\label{seq}
0 \lto \calm_{*} \otimes_{\calrs\otimes 1} \Omega^1_{\calrs/\calss}
\lto \calm_{*} \otimes_{\calrs\otimes 1} \calss / \cali^2_{*}
\xto{\id\otimes\mu}
\calm_{*} \lto 0 .
\end{equation}
that thus defines an extension class
$$
\ex_{\calm_{*}} \in \Ext^1_{\calrs}
\big(\calm_{*} , \calm_{*} \otimes_{\calrs}
\Omega^1_{\calrs} \big)\,.
$$

\begin{prop}\label{4.2.7}
$\ex_{\calm_{*}}=\At_{\calm_{*}}$.
\end{prop}

\begin{proof}
In view of \ref{4.2.5}(\ref{4.2.5.4}) we may assume that there is 
a connection $\nabla : \calm_{*} \to \calm_{*} \otimes_\calrs
\Omega^1_{\calrs / \caloy}$ on $\calms$.  The morphism
\begin{eqnarray*}
\xi : \calm_{*} &\longrightarrow &\calm_{*} \otimes_{\calrs}
\calss / \cali^2_{*} \\
m &\lto & m \otimes 1 + \nabla(m)\ .
\end{eqnarray*}
is easily seen to be an $1\otimes\calrs$-linear right-inverse to
$\id\otimes\mu$ in (\ref{seq}). Thus, the exact sequence (\ref{seq}) 
of DG modules is semi-split, and its extension class in
$\Ext^1_{\calrs}(
\calm_{*} , \calm_{*} \otimes \Omega^1_{\calrs/Y})$ is given by the
morphism
$$
[\partial , \xi] =
[\partial , \id_{\calm_{*}} \otimes 1 ] +
[\partial , \nabla ] = [\partial , \nabla ]:
\calm_{*} \lto \calm_{*} \otimes \Omega^1_{\calrs/Y}\,,
$$
see  \cite[X.126, Cor.1(b)]{Alg}.
As $[\partial , \nabla ]$ represents the Atiyah class, the result
follows.
\end{proof}

\subsection{Atiyah classes on $D(X)$}
\label{ATX}
We will now descend from $D(\calrs)$ to $D(X)$ to construct Atiyah
classes for any complex in $D(X)$.

Let $X$ be a complex space and let $(X_{*},W_{*},\calrs,\calbs)$ be an
extended resolvent of $X$, that is, we consider first the absolute case, where $Y$
is reduced to a point. With
$\bar\calbs:=\calbs\otimes_{1\otimes\calrs}\caloxs$, we showed in 
\cite[3.3.3.(3)]{BFl2} that 
for any complex of $\calox$--modules $\calm$, with $\calms$ the induced 
complex of $\caloxs$--modules on the simplicial space $X_{*}$, 
the morphism of complexes of right $\caloxs$--modules
$$
\calms\otimes_{\calrs\otimes 1} \bar\calbs
\xto{\id_{\calms}\otimes \bar\nu} \calms
$$
is a quasiisomorphism. Moreover,
$\Omega^1_{\calrs/Y}\otimes_\calrs \caloxs\cong\bbbl_{X_{*}}$ represents
the cotangent complex of $X_{*}$. Thus, applying $-\otimes_{1\otimes\calrs}\caloxs$
to $\id_\calms\otimes\At_\calbs$, as defined in \ref{4.2.4}, results in a morphism 
$\At_\calms$ defined by the diagram
\begin{align}
\label{eq11}
\bdi
\calms\otimes_{\calrs\otimes 1}\bar\calbs
&\rTo^{\id_{\calms}\otimes\At_\calbs\otimes 1} &
  \calms\otimes_{\calrs\otimes
1}\bar\calbs\otimes_{\caloxs}\bbbl_{X_{*}}[1]\\
\dTo^{\id_{\calms}\otimes \bar\nu}_{\simeq}&& 
\dTo_{\id_{\calms}\otimes \bar\nu\otimes \id}^{\simeq}\\
\calms & \rTo^{\At_{\calm_{*}}} &
\calms\otimes_{\caloxs} \bbbl_{X_{*}}[1]\,.
\edi
\end{align}
We will call $\At_\calms$ the {\em absolute Atiyah class\/} of $\calms$. 
Taking \Cech complexes and observing that $\calcb(\calms)\simeq \calm$ 
and $\calcb(\calms\otimes\bbbl_{X_{*}}[1])\simeq 
\calm\dotimes_\calox \bbbl_X[1]$, we obtain the {\em absolute Atiyah class\/} 
of $\calm$,
$$
\At_\calm\simeq \calcb(\At_\calms):\calm\to \calm\dotimes_{\calox} \bbbl_{X}[1]\,,
$$
as a morphism in $D(X)$. In the same way we can define the
{\em absolute Atiyah-Chern character\/} of $\calm$,
$$
\AC_\calm=\calcb(\exp(-\At_\calms))=:\exp(-\At_\calm): 
\calm\to \calm\dotimes_{\calox} \bbbs(\bbbl_{X}[1])\,.
$$

Now we turn to the relative case of an arbitrary morphism $f:X\to Y$
of complex spaces. Note that there is a natural comparison morphism 
$\bbbl_{X} \lto  \bbbl_{X/Y}$ from the absolute cotangent complex to 
the relative one. 

\begin{defn}\label{4.3.1}
For a morphism $X\to Y$ of complex spaces and $\calm\in D(X)$, 
we call the composite morphism 
$$
\At^{X/Y}_\calm:\calm\lto \calm\dotimes_{\calox} \bbbl_{X}[1] \lto
\calm\dotimes_{\calox} \bbbl_{X/Y}[1]
$$
in $D(X)$ the {\em Atiyah class} of $\calm$ {\em relative to\/} $X/Y$, and 
$$
\AC^{X/Y}_\calm:= \exp(-\At^{X/Y}_\calm): 
\calm\lto \calm\dotimes_{\calox}\bbbs( \bbbl_{X}[1]) \lto
\calm\dotimes_{\calox} \bbbs(\bbbl_{X/Y}[1])
$$
the {\em Atiyah-Chern character\/} of $\calm$ relative to $X/Y$.
\end{defn}

As before we will almost always  omit the upper indices of the Atiyah class and 
Atiyah-Chern character and write simply $\At_\calm$ and $\AC_\calm$. 

For complexes in $D^-_c(X)$, this
definition coincides with the one given in
\cite[Sect.\ 3]{BFl}, by \ref{4.2.5}(2). Finally, we note that 
\ref{At transformation} descends as well to $D(X)$, that is,
we have the following.

\begin{prop}\label{4.3.2}
For every morphism $f:X\to Y$, forming the Atiyah classes $\At_{\calm}$
or the Atyah-Chern character $\AC_\calm$ defines exact transformations
\begin{align*}
\At_{?}&:\id_{D(X)}\lto \id_{D(X)}\dotimes_{\calox}\bbbl_{X/Y}[1]
\quad\text{and}\\
\AC_{?}&:\id_{D(X)}\lto 
\id_{D(X)}\dotimes_{\calox}\bbbs_{\calox}(\bbbl_{X/Y}[1])
\end{align*}
of endofunctors on $D(X)$. In detail,

\begin{enumerate}[\quad\rm (1)]
\item
 Every morphism $\vp:\calm\to \caln$ in $D(X)$ fits into commutative
diagrams\footnote{The diagram on the left is, of course, just a piece of the diagram 
on the right.}
\begin{align*}
\bdi
\calm & \rTo^{\At_\calm} &\calm\dotimes_\calox \bbbl_{X/Y}[1]\\
\dTo^{\vp} && \dTo_{\vp\otimes\id}\\
\caln & \rTo^{\At_\caln}  &\caln\dotimes_\calox \bbbl_{X/Y}[1].
\edi
\qquad\text{and}\qquad
\bdi
\calm & \rTo^{\AC_\calm} &\calm\dotimes_\calox\bbbs_{\calox}(\bbbl_{X/Y}[1])\\
\dTo^{\vp} && \dTo_{\vp\otimes\id}\\
\caln & \rTo^{\AC_\caln} &\caln\dotimes_\calox \bbbs_{\calox}(\bbbl_{X/Y}[1])\,.
\edi
\end{align*}

\item
The functor $\At_{?}$ anticommutes with the translation functor, that is,
$\At_{\calm[1]}=-(\At_\calm)[1]$ for every complex $\calm\in D(X)$.
Equivalently, for every distinguished triangle 
$$
\Delta=(\calm'\xto{u}\calm\xto{v}\calm''\xto{w}\calm[1])
$$ 
in $D(X)$,
\begin{align*}
(\At_{\calm'},\At_{\calm},\At_{\calm''})&:
\Delta\lto \Delta\dotimes_\calox \bbbl_{X/Y}[1]\quad\text{and}\\
(\AC_{\calm'},\AC_{\calm},\AC_{\calm''})&:
\Delta\lto \Delta\dotimes_\calox \bbbs_{\calox}(\bbbl_{X/Y}[1])
\end{align*}
constitute morphisms of distinguished triangles in $D(X)$.
\end{enumerate}
\end{prop}

\bproof
Commutativity of the diagram (\ref{eq11}) above implies immediately 
that forming the Atiyah class or Atiyah-Chern character of complexes 
on $X_{*}$ is natural with respect to morphisms $\calms\to\calns$ in $D(X_{*})$,
and anticommutation with the translation functor descends from 
\ref{4.2.5}(\ref{4.2.5.4}) in $D(\calrs)$ to $D(X_{*})$. 

Applying the \Cech construction, the same properties are passed on 
to the absolute Atiyah class
$\At_\calm:\calm\to\calm\dotimes_\calox\bbbl_X$ in $D(X)$. As the (relative)
Atiyah class $\At_\calm:\calm\to\calm\dotimes_\calox\bbbl_{X/Y}$ is
obtained from the absolute one by composing with the
natural map
$\calm\dotimes_\calox\bbbl_{X} \to \calm\dotimes_\calox\bbbl_{X/Y}$,
the result follows in full generality.
\eproof

\section{The  Atiyah-Hochschild Character and the Centre
of the Derived Category}

\subsection{The Atiyah-Hochschild character}

In this subsection we will define the Atiyah-Hochschild character that assigns to every 
$\calm\in D(X)$ a canonical morphism
$$
\AH_\calm: \calm \lto \calm \dotimes_\calox\bbbh_{X/Y}
$$
in $D(X)$, and we will relate it to the Atiyah-Chern character $\AC_\calm$.

\bsit \label{AHdef}
Since we have a natural map $\bbbh_{X}\to \bbbh_{X/Y}$ it suffices to construct this 
character in case $Y$ is reduced to a point. 
As in \ref{ATX}, let us choose a free extended resolvent $\fX = (X_{*},W_{*},\calrs,\calbs)$
that comes with a quasiisomorphism $\nu:\calb\to\calrs$ of DG $\calss$--modules.
Using the \Cech functor (see \cite[2.27]{BFl}) it suffices further to construct such 
Atiyah-Hochschild characters 
$$
\AH_\calms:\calms \lto \calms \dotimes_\caloxs\bbbh_{X_*} 
$$
for $\calms\in D(X_*)$ with the usual functoriality properties. 

As in ({\em loc.cit.}) we set $\bar\calbs:= \calbs\otimes_{1\otimes\calrs}\caloxs$ and 
denote $\bar\nu: \bar\calbs\to\caloxs$ the ca\-no\-nical projection induced by $\nu$. 
The DG algebra homomorphism $1\otimes\calrs\to \calss$ induces a canonical map 
$$
q:\bar \calbs\to \bbbh_{X_*}\cong \calbs\otimes_{\calss} \caloxs\,,
$$
and the latter is a morphism of complexes of locally free $\caloxs$-modules. 
According to \cite[Lemma 3.3.3]{BFl2}, the map 
$\id\cdot\bar\nu: \calms\otimes_{\calrs\otimes 1}\bar\calbs\to \calms$ is a 
quasiisomorphism. We now define $\AH_\calms$ as a morphism in the derived 
category $D(X_*)$ through the diagram
\begin{align}
\label{eq61}
\bdi
\calms\otimes_{\calrs\otimes 1}\bar\calbs&\rTo^{\id\otimes q}&
\calms\otimes_{\caloxs}(\calbs\otimes_{\calss}\caloxs)\\
\dTo^{\id\cdot  \bar\nu}_{\simeq} &&
\dTo_{\simeq}\\
\calms  &\rTo^{\AH_\calms}&\calms\otimes_\caloxs \bbbh_{X_*}
\edi
\end{align}
It is easy to see from this construction that the following hold.
\esit

\bpro\label{AHpro1}
Assigning to an $\calox$-module the Atiyah-Hochschild character  $\AH_\calm$ is 
functorial in $\calm$ and anticommutes with the shift functor or, equivalently, is 
compatible with long exact Ext-sequences. 
In other words, each morphism $X\to Y$ of complex spaces defines an exact morphism
\begin{align*}
\AH_{?}:\id_{D(X)}\lto \id_{D(X)}\dotimes_{\calox}\bbbh_{X/Y}
\end{align*}
of endofunctors on $D(X)$.
\qed
\epro

For a morphism $f:X\to Y$ of complex spaces or schemes over a field of 
characteristic zero, the decomposition theorem yields via \ref{main4} a natural 
quasiisomorphism
$$
\Phi_{X/Y}:\bbbh_{X/Y}\xto{\ \simeq\ }\bbbs(\bbbl_{X/Y}[1])\,.
$$
Comparing via this map the Atiyah-Hochschild character and the Atiyah-Chern 
character we get the following result.

\bthm\label{AHthm}
Via the decomposition theorem the Atiyah-Hochschild character identifies with 
the Atiyah Chern character, that is, for any complex of $\calox$-modules $\calm$ 
the diagram 
\bdi
&& \calm  \\
&\ldTo^{\AH_\calm}&&\rdTo^{\AC_\calm}\\
\calm\dotimes_\calox \bbbh_{X/Y} && \rTo^{\id\otimes \Phi_{X/Y}}&&
\calm\dotimes_{\calox}\bbbs(\bbbl_{X/Y}[1])
\edi
commutes in $D(X)$.
\ethm

\bproof
In view of the definitions and \ref{4.3.2}, it suffices again to establish the absolute 
case, when $Y$ is reduced to a point. By definition of $\Phi_\calbs$; see \ref{stage}; 
the diagram
\bdi[midshaft]
\calbs&\rTo^{\AC_\calbs}&
\calbs\otimes_{\calss} \bbbs(\Omega^1_{\calss/\calrs \otimes 1}
[1])&\rEqual & \calbs\otimes_{1\otimes \calrs} \bbbs(\Omega^1_{\calrs}
[1])\\
\dTo<{\id\otimes 1} && && \dTo>{\nu\cdot \id}\\
\calbs\otimes_{\calss}\calrs&& \rTo^{\Phi_\calbs}&&
\bbbs(\Omega^1_{\calrs}[1])
\edi
commutes. Tensoring with 
$\calms\otimes_{\calrs\otimes 1}- \otimes_{1\otimes\calrs}\caloxs$ and using the 
identifications 
$$
\bbbs(\Omega^1_{\calrs}
[1])\otimes_\calrs\caloxs\cong \bbbs(\bbbl_{X_*}[1])\quad\mbox{and}\quad 
\calbs\otimes_{\calss}\caloxs\cong \bbbh_{X_*}
$$
we obtain that the outer square in the diagram
\bdi[midshaft]
 \calms\otimes_{\calrs\otimes 1}\bar \calbs
&&\rTo^{\id\otimes\AC_\calbs\otimes 1}  && 
\calms\otimes_{\calrs\otimes 1}\bar\calbs\otimes_{\caloxs} 
\bbbs(\bbbl_{X_*}[1])\\
&\rdTo^{\id\cdot \bar\nu}\\
\dTo<{\id\otimes q} && \calms&& 
\dTo>{\id\otimes (\bar\nu\cdot \id)}\\
&\ldTo^{\AH_\calms}&&\rdTo^{\AC_\calms}&\\
\calms\otimes_{\caloxs}\bbbh_{X_*}&& \rTo^{\id\otimes \Phi_{X_*}}&&
\calms\otimes_\caloxs\bbbs(\bbbl_{X_*}[1])
\edi
is commutative. By definition of $\AH_\calms$ and $\AC_\calms$ the triangle on the 
left and the right upper part of the diagram are commutative as well. 
As $\id\cdot\bar\nu$ is a quasiisomorphism, the result follows. 
\eproof

\subsection{The characteristic homomorphism to the graded centre}
In this part we will relate the Atiyah-Hochschild character to the mapping from 
Hochschild cohomology into the graded centre of the derived category.

\begin{sit}
\label{char}
Recall from \cite[Sect.3]{BFl2} that for every morphism $f:X\to Y$ of analytic 
spaces or schemes, not necessarily in characteristic zero, there is a natural 
homomorphism of graded commutative algebras from the Hochschild cohomology 
ring $\HH^{\bdot}_{X/Y}(\calox)$ to the graded centre of the derived category of $X$, 
a morphism we denoted
$$
\chi=\chi_{X/Y}: \HH^{\bdot}_{X/Y}(\calox)\lto \fZ_{gr}^{\bdot}(D(X))\,.
$$
Recall as well that in case $X$ flat over $Y$, this morphism is easily defined by 
means of the (identical) Fourier-Mukai transformation that is induced by the 
diagonal embedding $X\subseteq X\times X$.

Each complex $\calm\in D(X)$ gives rise to an evaluation map $ev_{\calm}:
\fZ_{gr}^{\bdot}(D(X))\to \Ext_X^{\bdot}(\calm,\calm)$, a homomorphism of graded rings 
that takes its values in the graded centre, thus, endows $\Ext_X^{\bdot}(\calm,\calm)$
with the structure of a graded $\fZ_{gr}^{\bdot}(D(X))$--algebra.

To give the map $\chi$ amounts then to describe the family 
$(\chi_{\calm}=ev_{\calm}\circ \chi)_{\calm}$ 
of homomorphisms of graded rings
$$
\chi_{\calm}:\HH^{\bdot}_{X/Y}(\calox) \lto \Ext_X^{\bdot}(\calm,\calm)\,,
$$
natural both with respect to morphisms in $D(X)$ and to distinguished triangles,
equivalently, with respect to the long exact $\Ext$--sequences induced by such triangles.
\end{sit}

Contracting against the Atiyah-Hochschild character of a complex $\calm$
yields as well a map
$$
?\lr\AH_\calm:\HH^\bdot_{X/Y}(\calox)\cong \Ext^{\bdot}_X(\bbbh_{X/Y},\calox)
\lto \Ext_X^{\bdot}(\calm,\calm)
$$
that assigns to a morphism $\beta:\bbbh_{X/Y}\to\calox[r]$ in D(X) the composition
$$
\beta\lr\AH_\calm :
\calm\xto{\AH_\calm}\calm\dotimes_{\calox}\bbbh_{X/Y}
\xto{\id_{\calm}\otimes\beta}\calm\dotimes_{\calox}\calox[r]\cong \calm[r]\,.
$$

\bpro
\label{module}
For every complex of $\calox$-modules $\calm$ the map $\chi_\calm$ is equal to the 
contraction against the Atiyah-Hochschild character
$$
\chi_\calm=(?\lr \AH_\calm): \HH^{\bdot}_{X/Y}(\calox) \lto \Ext_X^{\bdot}(\calm,\calm)\,.
$$
\epro

\bproof
Let $\beta\in \HH^r_{X/Y}(\calox)$ be represented by a morphism 
$\beta: \bbbh_{X/Y}\to\calox$ in $D(X)$. We need to show that the composed morphism
$$
\calm\xto{AH_\calm} \calm\dotimes_\calox \bbbh_{X/Y}\xto{\id\cdot \beta}
\calm
$$
is equal to $\chi_\calm.$
As before we may assume that $Y$ is reduced to a point. 
To achieve the aim, we use the notations of \ref{AHdef}. 
It is sufficient to show that $\chi_\calms= (?\lr\AH_\calms)$ on $X_*$. 

As $\calbs$ is generally not a projective $\calss$--module, we will choose first a projective 
approximation $\calps\to\calbs$. Tensoring this morphism of flat DG $\calss$--modules 
with $\caloxs$ returns a quasiisomorphism 
$\bar\calps:=\calps\otimes_{1\otimes\calrs}\caloxs\to\bar\calbs$, and composition with 
$\bar\nu$ leads to a quasiisomorphism, say ${\bar\mu}:\bar\calps\to\caloxs$. 
As $\id\cdot\bar\nu:\calms\otimes_{\calrs\otimes 1}\bar\calbs\to\calms$ is  a 
quasiisomorphism, the same holds for 
$\id\cdot{\bar\mu}:\calms\otimes_{\calrs\otimes 1}\bar\calps\to\calms$. 

The element $\beta$ is represented by an actual  morphism 
$\bar\beta:\bar\calps\to\bar\calps[r]$ of complexes of $\caloxs$--modules under the 
identification $\HH^r_{X}(\calox)\cong H^0(\Hom_\caloxs(\bar\calps,\bar\calps[r]))$. 
Now consider the commutative diagram
\bdi[midshaft]
\calms&\lTo^{\id\cdot{\bar\mu} }&\calms\otimes_\calrs\bar\calps &\rTo^{p=proj}& 
\calms\otimes_\calrs\calps\otimes_{\calss} \caloxs&\rEqual & 
\calm\dotimes_\caloxs \bbbh_{X_*}\\
\dTo>{\chi_\calms(\beta)} && \dTo>{\id\otimes\bar\beta} &&
\dTo>{\id\cdot \beta}\\
\calms[r]&\lTo^{\id\cdot{\bar\mu}}&\calms\otimes_\calrs\bar\calps [r] &
\rTo^{\id\cdot{\bar\mu}}& \calms[r]\,.
\edi
As a morphism in the derived category, $\chi_\calms(\beta)$ is defined by the 
square on the left because $\id\cdot{\bar\mu}$ is a quasiisomorphism. 
By definition, $\AH_\calms$ equals $p\circ(\id\cdot{\bar\mu})^{-1}$ in the derived 
category. Thus, the result follows from the commutativity of the diagram.
\eproof

For a morphism $f:X\to Y$ of complex spaces or schemes over a field of 
characteristic zero, the decomposition theorem yields via \ref{main4} a 
natural isomorphism
$$
\Psi := \Hom_{D(X)}(\Phi_{X/Y},\calox):\Ext^{\bdot}_X(\bbbs(\bbbl_{X/Y}[1]),\calox)
\xto{\ \cong\ } \HH^{\bdot}_{X/Y}(\calox)\,,
$$
while contracting against the Atiyah-Chern character of a complex $\calm$
yields a map
$$
?\lr\AC_{\calm}:\Ext^{\bdot}_X(\bbbs(\bbbl_{X/Y}[1]),\calox)
\lto \Ext_X^{\bdot}(\calm,\calm)
$$
Using \ref{AHthm} and \ref{module} we obtain immediately the following result.

\begin{cor}\label{semiregpro1}
The diagram 
\bdi[height=10mm]
\Ext^{\bdot}_X(\bbbs(\bbbl_{X/Y}[1]),\calox)&&\rTo^{\Psi}_\cong&& \HH^{\bdot}_{X/Y}(\caloxs)\\
&\rdTo_{?\lr \AC_\calm}&&\ldTo_{\chi_\calm}\\
& &\Ext_X^\bdot(\calm,\calm)
\edi
commutes.\qed
\end{cor}

\subsection{The map from Hochschild cohomology into 
$\Ext^{\bdot}_{X}(\bbbh_{X/Y},\bbbh_{X/Y})$} 

\bdfn\label{5.3.1}
In connection with the Atiyah-Hochschild character it is also useful to study the 
canonical ring homomorphism  
$$
\alpha: \HH^\bdot_{X/Y}(\calox)\to \Ext^\bdot_X(\bbbh_{X/Y},\bbbh_{X/Y})
$$
that we introduce as follows. With $\fX = (X_{*},W_{*},\calrs,\calbs)$ a free 
extended resolvent we have canonical isomorphisms
$$\ba{c}
\HH^\bdot_{X/Y}(\calox)\cong \Ext^\bdot_{\calss}(\calbs,\calbs)\,,
\\
\Ext^\bdot_X(\bbbh_{X/Y},\bbbh_{X/Y})\cong 
\Ext^\bdot_\calrs(\calbs\otimes_{\calss}\calrs, \calbs\otimes_{\calss}\calrs).
\ea
$$
Now we define $\alpha$ through the canonical homomorphism of Ext--algebras 
$$
\Ext^\bdot_{\calss}(\calbs,\calbs)\to
\Ext^\bdot_\calrs(\calbs\otimes_{\calss}\calrs,\calbs\otimes_{\calss}\calrs)
$$
induced by base change along the multiplication map $\calss\to\calrs$.
\edfn

\brem
\label{5.3.2}
(1) The multiplication on 
$\HH^\bdot_{X/Y}(\calox)\cong  \Ext^\bdot_X(\bbbh_{X/Y},\calox)$ 
can be obtained as the composition of the maps 
$$\ba{lll}
\HH^\bdot_{X/Y}(\calox)\times \HH^\bdot_{X/Y}(\calox)
&\xto{can \times \alpha}&
 \Ext^\bdot_X(\bbbh_{X/Y},\calox)\times \Ext^\bdot_X(\bbbh_{X/Y},\bbbh_{X/Y})\\
&\xto{\phantom{\alpha \times can}}& 
\Ext^\bdot_X(\bbbh_{X/Y},\calox)\cong \HH^\bdot_{X/Y}(\calox)\,,
\ea
$$
where the second arrow is the Yoneda product. 
Indeed, using as in definition \ref{5.3.1} the isomorphism 
$\HH^\bdot_{X/Y}(\calox)\cong \Ext^\bdot_{\calss}(\calbs,\calbs)$ this follows 
easily from the fact that the product is given by composition of morphisms.

(2)  Similarly, the left $\HH^\bdot_{X/Y}(\calox)$--module structure on Hochschild 
homology $\HH_{-{\scriptstyle \star}}^{X/Y}(\calox) \cong 
\Ext^{\scriptstyle \star}_X(\calox,\bbbh_{X/Y})$ is the composition  of the maps 
$$\ba{lll}
 \HH^\bdot_{X/Y}(\calox)\times\HH_{-{\scriptstyle \star}}^{X/Y}(\calox)
&\xto{\alpha \times can }&
\Ext^\bdot_X(\bbbh_{X/Y},\bbbh_{X/Y})\times 
\Ext^{\scriptstyle \star}_X(\calox,\bbbh_{X/Y})\\
&\xto{\phantom{can \times\alpha}}& 
\Ext^{\bdot+{\scriptstyle \star}}_X(\calox, \bbbh_{X/Y})\cong 
\HH_{-\bdot-{\scriptstyle \star}}^{X/Y}(\calox)\,,
\ea
$$
where again the second arrow is the Yoneda product. The argument is 
similar as in (1) and is left to the reader.
\erem

Recall that the natural ring homomorphism 
$\chi:\HH^{\bdot}_{X/Y}(\calox)\to \fZ_{gr}^{\bdot}(D(X))$ into the centre 
of the derived category  induces for every complex $\calm$ in $D(X)$ a 
natural ring homomorphism
$$
\chi_\calm:\HH^{\bdot}_{X/Y}(\calox)\lto \Ext_X^\bdot(\calm, \calm)\,,
$$
see \ref{module}.
In particular, $\Ext_X^\bdot(\calm, \calm)$ is endowed with  a module 
structure over  $\HH^{\bdot}_{X/Y}(\calox)$ that we denote $f.b$ for 
$b\in\HH^{\bdot}_{X/Y}(\calox)$ and $f\in \Ext_X^\bdot(\calm, \calm)$. 
Moreover, Hochschild homology $\HH^{X/Y}_{\bdot}(\calox)$ also carries 
a natural structure of a $\HH^{\bdot}_{X/Y}(\calox)$-module. We next 
establish that the Atiyah-Hochschild character is compatible with these structures.

\bpro\label{5.3.3}
For every complex $\calm\in D(X)$ and for any elements 
$b\in\HH^{\bdot}_{X/Y}(\calox)$ and $f\in \Ext^\bdot_X(\calm,\calm)$ 
the diagram
\bdi
\calm  &\rTo^{\AH_\calm} & \calm\otimes_\calox \bbbh_{X/Y}\\
\dTo<{f.b} && \dTo>{f\otimes\alpha(b)} \\
\calm  &\rTo^{\AH_\calm} & \calm\otimes_\calox \bbbh_{X/Y}
\edi
commutes, thus, the Atiyah-Hochschild character is a homomorphism 
of (right) modules over the Hochschild cohomology ring.
\epro 

\bproof
As before we may assume that $Y$ is reduced to a point. 
To achieve our aim, we place ourselves in the context of \ref{ATX}, 
thus, work with a fixed free extended resolvent $\fX = (X_{*},W_{*},\calrs,\calbs)$
that comes with a quasiisomorphism $\nu:\calbs\to\calrs$ of DG $\calss$--modules. 
It is sufficient to deduce the analogous commutative diagram above on $X_*$. 
The map $\nu$ induces a quasiisomorphism 
$\bar\nu:\bar\calbs:=\calbs\otimes_{1\otimes\calrs}\caloxs\to\caloxs$, 
see \cite[Lemma 3.3.3]{BFl2}. 

As $\calbs$ is not a projective $\calss$--module, we choose again first a 
projective approximation $\calps\to\calbs$ as in \ref{module}. As there, 
tensoring with $\caloxs$ results in a quasiismorphism 
$\bar\calps:=\calps\otimes_{1\otimes\calrs}\caloxs\to\bar\calbs$, and 
composition with $\bar\nu$ leads to a quasiisomorphism 
${\bar\mu}:\calps\to\caloxs$.
Now the element $b$ is represented by a morphism 
$\bar b:\bar\calps\to\bar\calps[r]$ under the identification 
$\HH^r_{X}(\calox)\cong H^0(\Hom_\caloxs(\bar\calps,\bar\calps[r]))$. 
Next consider the  commutative diagram
\bdi
\calms&\lTo^{\id\cdot {\bar\mu}}&\calms\otimes_\calrs \bar\calps  &
\rTo^{\id\otimes proj}& \calms\dotimes_\caloxs \bbbh_{X_*}\\
\dTo>{f.b} && \dTo>{f\otimes \bar b} &&\dTo>{f\otimes\alpha(b)}\\
\calms&\lTo^{\id\cdot{\bar\mu}}&\calms\otimes_\calrs\bar\calps &
\rTo^{\id\otimes  proj}& \calms\otimes_\caloxs\bbbh_{X/Y}\,.
\edi
The morphism $f.b$ is in the derived category $D(X_*)$ the morphism 
defined through the square on the left. By definition, the compositions 
$(\id\otimes  proj)\circ (\id\cdot {\bar\mu})^{-1}$ of the maps along the top 
or bottom are both equal to $\AH_\calms$.  Using the definition of 
$\alpha(b)$ as in \ref{5.3.1} the result follows.
\eproof

\section{Hochschild (Co-)Homology and the Semiregularity Map}
\subsection {The semiregularity map revisited}
Semiregularity is a concept that pertains to {\em perfect\/} complexes only,
as it involves taking traces of endomorphisms. We review the pertinent
definition and results from \cite{BFl}, but first simplify notation:
Let $f:X\to Y$ be a fixed morphism of complex spaces and denote
$\bbbl:=\bbbl_{X/Y}$ its cotangent complex, and 
$\bbbs:= \bbbs_{\calox}(\bbbl_{X/Y}[1])$ the derived symmetric algebra
over the shifted cotangent complex. Furthermore, all unadorned (derived) 
tensor products in this section are to be taken over $\calox$.

We remind the reader that the cohomology of the cotangent complex is also called
the {\em tangent\/} or {\em Andr\'e-Quillen cohomology\/} of $X/Y$, denoted
$T^{\bdot}_{X/Y}(\calm) := \Ext^{\bdot}_{X}(\bbbl,\calm)$ for any
$\calm\in D(X)$.

\begin{sit}
Given a {\em perfect\/} complex of $\calox$-modules $\calf$, 
we can first assign to it its Atiyah-Chern character, the exponential of the 
Atiyah class
$$
AC_\calf=\exp(-\At_{\calf})\in \Ext^{0}_{X}(\calf, 
\calf\dotimes \bbbs)\,,
$$
and then take the trace to obtain the {\em Chern character\/}
$$
\ch(\calf):= \Tr (\AC_\calf)\in H^{0}(X,\bbbs)
$$ 
\`a la Illusie \cite{Ill} in general, or as in O'Brian-Toledo-Tong 
\cite{OTT} in the smooth complex analytic case.
\end{sit}

\begin{sit}
\label{semiregreview}
In \cite{BFl}, we introduced for any perfect complex $\calf$ on $X$ its 
associated {\em semiregularity map\/}
$$
\sigma_{\calf}: \Ext^{\bdot}_{X}(\calf, \calf)\lto H^{\bdot}(X,\bbbs)
$$
given, by taking the trace after multiplying with the Atiyah-Chern character, as
$$
\sigma_{\calf}(\vp):= \Tr\left(\exp(-\At_{\calf})\circ \vp : \calf
\xto{\vp}\calf\xto{\AC_{\calf}} \calf\dotimes_{\calox}\bbbs\right)\,.
$$
We then showed, as a special case of \cite[4.2]{BFl}, that the following diagram
\begin{align}
\label{semiregtriangle}
\bdi
T^{\bdot-1}_{X/Y}(\calox) = \Ext^{\bdot-1}_{X}(\bbbl,\calm)
&&\rTo^{?\lr (-\At_{\calf})}&&\Ext^{\bdot}_{X}(\calf,\calf)\\
&\rdTo<{?\lr \ch(\calf)}&&\ldTo_{\sigma_{\calf}}\\
&&H^{\bdot}(X,\bbbs)
\edi
\end{align}
commutes, thus, relating the tangent cohomology of $X$ over $Y$ to
the semiregularity map and the Chern character of a perfect complex.
\end{sit}

Our aim in this section is now twofold. First, we define a semiregularity  map into 
Hochschild homology using the Atiyah-Hochschild character. Second, we extend 
the diagram above, replacing tangent cohomology by Hochschild cohomology.

\subsection{The semiregularitay map on Hochschild cohomology}

We now show that the Atiyah-Hoch\-schild character induces a semiregularity 
map to Hochschild homology. Let $\calf$ be a perfect complex on $X$ and define 
$$
\sigma^H_\calf:\Ext^\bdot_X(\calf, \calf)\lto \HH_{-\bdot}^{X/Y}(\calox)
$$
through the formula
$$
\sigma^H_\calf(f):=\Tr ((f\otimes \id)\circ \AH_\calf)\in \Ext^\bdot_X(\calox,\bbbh_{X/S})\,, 
$$
where $(f\otimes \id)\circ \AH_\calf$ is the composition 
$$
\calf\xto{\AH_\calm} \calf\otimes\bbbh_{X/Y}\xto{f\otimes \id}
\calf\otimes\bbbh_{X/Y}\,.
$$
Comparing this with the semiregularity map $\sigma_\calf$ as introduced in 
\cite{BFl} we obtain the following result.

\bpro
For every perfect complex $\calf$ of $\calox$-modules the diagram
\bdi
&&\Ext^\bdot_X(\calf, \calf) \\
&\ldTo^{\sigma^H_\calf} && \rdTo^{\sigma_\calf}\\
\HH_{-\bdot}^{X/Y}(\calox) &&\rTo^{H^{\bdot}(X,\Phi_{X/Y})}&&
\Ext_X^\bdot (\calox, \bbbs(\bbbl_{X/Y}[1]))
\edi
commutes, where $\sigma_\calf$ is the semiregularity map.
\epro

\bproof
By Theorem \ref{AHthm} the diagram 
\bdi
\calf&\rTo^{\AH_\calm} &\calf\otimes\bbbh_{X/Y}&\rTo^{f\otimes \id}&
\calf\otimes\bbbh_{X/Y}\\
\dTo^{\id} &&\dTo_{\id \otimes\Phi_{X/Y}}&&\dTo_{\id \otimes\Phi_{X/Y}}\\
\calf&\rTo^{\AC_\calm} &\calf\otimes\bbbs(\bbbl_{X/Y}[1])&\rTo^{f\otimes \id}&
\calf\otimes\bbbs(\bbbl_{X/Y}[1])
\edi
commutes. Taking traces now yields the result. 
\eproof

\begin{rem}
As $\Ext^\bdot_X(\calf, \calf)$ is concentrated in nonnegative (cohomological) 
degrees and $\sigma^H_\calf$ preserves degrees, the semiregularity map into 
Hochschild homology $\HH_{-\bdot}^{X/Y}(\calox)$ can only be nonzero in 
{\em nonpositive\/} (homological) degrees, whence this contraction takes its 
values only in the nonpositive components 
$\bigoplus_{n\le 0} \HH_{n}^{X/Y}(\calox)$ of Hochschild homology.
\end{rem}

\subsection{The semiregularity map is linear over the Hochschild cohomology} 

Recall that the natural ring homomorphism 
$\chi:\HH^{\bdot}_{X/Y}(\calox)\to \fZ_{gr}^{\bdot}(D(X))$ into the centre  of the 
derived category  induces for every complex $\calf$ in $D(X)$ a natural ring 
homomorphism 
$$
\chi_\calm:\HH^{\bdot}_{X/Y}(\calox)\lto \Ext_X^\bdot(\calf, \calf)\,,
$$
see \ref{module}.
In particular, $\Ext_X^\bdot(\calf, \calf)$ is endowed with  a right and left module 
structure over  $\HH^{\bdot}_{X/Y}(\calox)$. Let  $f.b$ and $b.f$ denote the scalar 
product for $b\in\HH^{\bdot}_{X/Y}(\calox)$ and $f\in \Ext_X^\bdot(\calf, \calf)$. 
Since $\chi$ is a map into the graded centre of the derived category we have
$$
b.f=(-1)^{|b||f|}f.b\,.\leqno (*)
$$
Moreover, Hochschild homology $\HH^{X/Y}_{\bdot}(\calox)$ also carries a 
natural structure of a left $\HH^{\bdot}_{X/Y}(\calox)$-module, as seen in \ref{5.3.2}(2). 
These structures are compatible as the next result shows.

\bpro\label{6.3.1}
For every perfect complex $\calf$ the semiregularity map 
$$
\sigma^H_\calf:\Ext^\bdot_X(\calf, \calf)\lto \HH^{X/Y}_{-\bdot}(\calox)
$$
is a linear map of left $\HH^{\bdot}_{X/Y}(\calox)$-modules.
\epro  

\bproof
Let $b\in \HH^r_{X/Y}(\calox)$ be given. We know from \ref{5.3.3} that 
for $f\in\Ext^\bdot_X(\calf,\calf) $
$$
\AH_\calf \circ (f.b)=
(f\otimes \alpha (b))\circ \AH_\calf=(-1)^{|f||b|}
(\id\otimes \alpha (b))\circ ( f\otimes\id)\circ \AH_\calf.
$$
Using $(*)$ this gives 
$$
\AH_\calf \circ (b.f)=
(\id\otimes \alpha (b))\circ ( f\otimes\id)\circ \AH_\calf.
$$
By the usual properties of traces, see \cite{Ill}, it follows that 
$$
\Tr(\AH_\calf \circ (b. f))=
\Tr ((\id\otimes \alpha (b))\circ ( f\otimes\id)\circ \AH_\calf)=
\alpha (b) \circ \Tr( \AH_\calf \circ f)
$$
and so $\sigma^H_\calf(b.f)=\alpha(b).\sigma^H_\calf(f)=b.\sigma^H_\calf(f)$, 
where the last equality follows from Remark \ref{5.3.2}(2).
\eproof

In analogy with the Chern character we can now introduce for every perfect 
complex on $X$ the {\em Chern-Hochschild character\/} to be 
$$
\ch^H_\calf:=\sigma^H_\calf(\id_\calf)\,.
$$
In terms of this character Proposition \ref{6.3.1} implies the following result.

\bcor\label{AHcom}
For every perfect complex $\calf$ of $\calox$--modules the diagram
\bdi
\HH^\bdot_{X/Y}(\calox)&& \rTo^{\chi_\calf} && \Ext^\bdot_X(\calf,\calf) \\
& \rdTo_{\cdot\ch^H_\calf} && \ldTo_{\sigma^H_\calf}\\
&&\HH_{-\bdot}^{X/Y}(\calox)
\edi
commutes.\qed
\ecor

\subsection{The semiregularity map and the Chern character}
We finish this paper with the promised generalization of the triangle in 
\ref{semiregreview}(\ref{semiregtriangle}), replacing tangent cohomology by 
Hochschild cohomology as the source of the relevant maps.
\begin{sit}
To formulate this result, note that 
$\Ext^{\bdot}_{X}(\bbbs,\calox)\cong \HH^{\bdot}_{X/Y}(\calox)$ acts naturally 
from the left on $H^{\bdot}(X,\bbbs)\cong \HH_{-\bdot}^{X/Y}(\calox)$ by 
composition, an action we denote by $?\lr ??$. In particular, 
contracting against the Chern character of a perfect complex yields a 
degree preserving map
$$
?\lr \ch_{\calf}: \Ext^{\bdot}_{X}(\bbbs,\calox)\lto H^{\bdot}(X,\bbbs)\,.
$$
\end{sit}

The result is now the following.

\bthm
\label{semiregthm}
For every perfect complex $\calf$ on $X$, the diagram
\bdi
\Ext^{\bdot}_{X}(\bbbs,\calox)&\rTo^{\Psi}_{\cong}
&\HH^{\bdot}_{X/Y}(\calox)&\rTo^{\chi}&\fZ^{\bdot}_{gr}(D(X))\\
\dTo^{?\lr \ch_{\calf}}&&\dTo_{\cdot \ch^H_\calf}&&\dTo_{ev_{\calf}}\\
H^{\bdot}(X,\bbbs)&\lTo^{H^{\bdot}(X,\Phi_{X/Y})}_{\cong}
&\HH_{-\bdot}^{X/Y}(\calox)&\lTo^{\sigma^{H}_{\calf}}&\Ext^{\bdot}_{X}(\calf,\calf)
\edi
commutes.
\ethm

\brem
By Corollary \ref{AHcom} the right-hand square in the diagram commutes. 
Therefore the additional information is that also the square on the left is 
commutative. Note that $\Phi_{X/Y}$, and then also $H^{\bdot}(X,\Phi_{X/Y})$, 
is {\em an isomorphism of rings\/}, but, in general, the isomorphism $\Psi$ will 
{\em not respect the multiplicative structure}; see \cite{BW} for explicit examples, 
where the decomposition theorem is not compatible with the ring structure. 
Thus, commutativity of the square on the left says, inter alia, that the difference 
in the ring structures cannot be detected under taking products with or contraction 
against (Hochschild-)Chern characters.
\erem

The {\em Proof of \ref{semiregthm}\/} will follow from Corollary \ref{semiregpro1} above
that describes the composition along the top and down the right-hand side 
in the diagram and from Proposition \ref{semiregpro2} below that we formulate in 
slightly greater generality than is required for the proof here.

The strategy of proof is the same as for our earlier result \cite[4.2]{BFl}, except that
we have to extend the considerations there from derivations to differential operators.

We thus set up first the framework as follows. The group
$$
A^{\bdot}:= \Ext^{\bdot}_X(\calf ,
\calf\dotimes\bbbs)
$$
carries a natural graded algebra structure that is associative, but in
general not graded commutative. Note further that the diagonal morphism
$\bbbl\to \bbbl\oplus \bbbl$ induces a (co-associative and cocommutative) 
comultiplication $\bbbs\to\bbbs\otimes\bbbs$.

An element $\xi\in \Ext^r_X(\bbbs^{k},\caln)$, for $\caln$ any $\calox$--module, 
defines now a {\em differential operator\/}
$$
\xi\lr ?: A\to \Ext^{\bdot}_X(\calf ,
\calf\dotimes\caln\dotimes\bbbs)
$$ 
of degree $r$ and order $k$ that is induced by the composition
\bdi
\bbbs &\rTo^{\Delta} & \bbbs^k\dotimes
\bbbs &\rTo^{\xi\otimes \id} & 
\caln[r]\dotimes\bbbs\,,
\edi
where $\Delta$ is the indicated component of the comultiplication on $\bbbs$.

The  divided $k$th power of the Atiyah class of $\calf$ on the other hand determines 
a morphism $(-1)^k\At^k_\calf/k!:\calf\to \calf\dotimes\bbbs^k$ of degree $0$ in the 
derived category. 
Contraction against this class results then in a morphism
\bdi[w=7mm,h=11mm]
\Ext_X^r(\bbbs^k,\caln)&& \rTo^{?\lr (-1)^k\At^k_\calf/k!} &&
\Ext_{X}^r(\calf, \calf\dotimes\caln)\,.
\edi
More concretely, a morphism $\xi\in \Ext_X^r(\bbbs^k,\caln)$ is mapped to the composition 
\bdi
\xi\lr ({(-1)^k\At^k\calf}/{k!}):\calf &\rTo^{(-1)^k\At^k_\calf/k!}& \calf \dotimes \bbbs^k & 
\rTo^{\id\otimes \xi} &\calf\dotimes \caln[r]\,,
\edi
when viewing elements of Ext-groups as morphisms in the derived category. 
In analogy with Proposition 4.2 in \cite{BFl}  we will now show the following result.

\begin{prop}  
\label{semiregpro2}
Let $X\to Y$ be a morphism of complex spaces. If $\calf$ is a perfect
complex on $X$ and $\caln$ is an $\calox$--module then the diagram
\bdi[w=7mm,h=11mm]
\Ext_{X}^{\bdot}(\bbbs^k,\caln)&& \rTo^{?\lr (-1)^k\At^k_\calf/k!} &&
\Ext_{X}^{\bdot}(\calf, \calf\dotimes\caln)\\
& \rdTo_{?\lr \ch(\calf)}&&
\ldTo^{\sigma_{\calf}^{\caln}}\\ 
&&H^{\bdot}(X,\caln\dotimes\bbbs)
\edi
commutes, where $\sigma_{\calf}^{\caln}$ denotes the semiregularity map, generalized
in a straightforward manner to accommodate the insertion of $\caln$ as in \cite[4.1]{BFl}.
\end{prop}

\begin{proof}
The argument is similar to that in the proof of Proposition 4.2 in \cite{BFl}.
We need to show that, for a fixed $\xi \in \Ext^r(\bbbs^k,\caln)$ and for each $n\ge 0$, 
$$
\Tr\left(\frac{\At^n_\calf}{n!}\cdot \left(\xi \lr \frac{\At_\calf^k}{k!}\right)\right)
= {\xi}\lr \Tr\left( \frac{\At^{k+n}_{\calf}}{(k+n)!}\right)\,.
$$
As the trace map is compatible with taking cup
products, for any integer $s$, the diagram
\bdi[h=7mm]
\Ext_{X}^{s}(\calf,  \calf\dotimes
\bbbs^{k+n}) &\rTo^{\Delta} &
\Ext_{X}^{s}(\calf, \calf\dotimes\bbbs^k\dotimes \bbbs^n)
&\rTo^{\xi} &
\Ext_{X}^{s+r}(\calf, \calf\dotimes\caln\dotimes\bbbs^n
)\\
\dTo>\Tr && \dTo>\Tr  && \dTo>\Tr \\
H^{s}(X,\bbbs^{k+n})
&\rTo^{\Delta} &
H^{s}(X,\bbbs^k\dotimes \bbbs^n)
&\rTo^\xi &
H^{s+r}(X, \caln\dotimes\bbbs^n)
\edi
commutes.  Therefore
\begin{align*}
{\xi}\lr \Tr\left( \frac{\At^{k+n}(\calf)}{(k+n)!}\right) =  
\Tr\left({\xi}\lr \frac{\At^{k+n}(\calf)}{(k+n)!}\right)
\end{align*}
Thus it suffices to show, for each $m\ge 0$,  that
\begin{align*}\tag{1}
\Tr\left(\frac{\At^{n+m}_\calf}{n!}\cdot 
\left(\xi \lr \frac{\At_\calf^k}{k!}\right)\right)=
\Tr\left({\At^{m}_\calf}\cdot 
\left(\xi \lr \frac{\At_\calf^{k+n}}{(k+n)!}\right)\right)\,.
\end{align*}
As $\xi\lr ?$ is a differential operator of degree $r$ and order $k$ on
$A^{\bdot}=\Ext^{\bdot}_{X}(\calf, \calf\dotimes \bbbs)$ with values in 
$M^{\bdot}:=\Ext^{\bdot}_{X}(\calf, \calf\dotimes\caln \dotimes\bbbs)$, 
and as $\At^k_\calf\in A^0$, we obtain that with $x:=\At_\calf$ the 
map $\xi':A^{\bdot}\to M^{\bdot} $ given by 
$$
A^{\bdot}\ni \eta\mapsto \xi\lr (x\eta) -
x(\xi\lr \eta)\in M^{\bdot}
$$
is a differential operator in $\eta$ of order $k-1$. The case $k=0$ being trivial,
we may use induction on $k\ge 0$ and assume that we have already shown
$$
\Tr\left(\frac{x^{n+m}}{n!}\cdot\left(\xi'\lr \frac{x^{k-1}}{(k-1)!}\right)\right)
= \Tr \left( x^m\cdot \left(\xi' \lr \frac{x^{k+n-1}}{
(k+n-1)!}\right) \right)
\leqno (2)
$$
The term on the left is by definition equal to
$$
\Tr\left(\frac{x^{n+m}}{n!}\cdot\left(\xi\lr \frac{x^{k}}{(k-1)!}\right)\right) 
\leqno (3)
$$
as $\xi\lr ( x^{k-1}/(k-1)!)=0$. 
Moreover, the term on the right in (2) is by definition equal to 
$$
\ba{l}
\Tr \left(x^m\cdot\left(\xi\lr  \frac{x^{k+n}}{
(k+n-1)!}\right)\right) -
\Tr \left(x^{m+1}\cdot \left( \xi\lr  \frac{x^{k+n-1}}{
(k+n-1)!}\right)\right) 
\ea
\leqno (4)
$$
Thus we obtain from (2), (3) and (4) the equation 
\begin{align*}
\tag{5}
&\Tr \left(x^m\cdot\left( \xi\lr  \frac{x^{k+n}}{
(k+n-1)!}\right)\right)\\
&=
\Tr \left( x^{m+1}\cdot \left( \xi\lr \frac{x^{k+n-1}}{
(k+n-1)!}\right)\right)
+\Tr\left(\frac{x^{n+m}}{n!}\cdot\left(\xi\lr \frac{x^{k}}{(k-1)!}\right)\right)\,.
\end{align*}
Using induction with respect to $n$, we may assume that 
$$
\Tr \left( x^{m+1}\cdot\left( \xi\lr  \frac{x^{k+n-1}}{
(k+n-1)!}\right)\right)=
\Tr \left(\frac{x^{n+m}}{(n-1)!} \left( \xi\lr  \frac{x^{k}}{
k!}\right)\right)\,.
$$
Comparing with (5), the claim (1) follows.
\end{proof}

\end{document}